\documentclass{amsart}
\usepackage{amssymb}
\usepackage{amsmath}
\usepackage{mathrsfs}
\usepackage{pb-diagram,lamsarrow,pb-lams,amscd}
\newenvironment{remark}{\vspace{4pt}\noindent\textit{Remark.}} {\qed \\[2mm]}
\newtheorem{theorem}{Theorem}[section]
\newtheorem{corollary}[theorem]{Corollary}
\newtheorem{prop}[theorem]{Proposition}
\newtheorem{lemma}[theorem]{Lemma}
\newtheorem{cor}[theorem]{Corollary}
\newtheorem{definition}[theorem]{Definition}

\newtheorem{Main theorem}[theorem]{Main theorem}

\newenvironment{example}[1]{\noindent\textbf{Example:} \textit{#1}}{\qed}
\newcommand{\LL}{\ensuremath{\mathbb{L}}}

\newcommand{\1}{\ensuremath{\mathbf{1}}}
\newcommand{\N}{\ensuremath{\mathbb{N}}}
\newcommand{\Z}{\ensuremath{\mathbb{Z}}}
\newcommand{\Q}{\ensuremath{\mathbb{Q}}}
\newcommand{\R}{\ensuremath{\mathbb{R}}}
\newcommand{\C}{\ensuremath{\mathbb{C}}}
\newcommand{\Pro}{\ensuremath{\mathbb{P}}}
\newcommand{\A}{\ensuremath{\mathbb{A}}}
\newcommand{\Ar}{\ensuremath{\mathcal{A}r}}
\newcommand{\Aa}{\ensuremath{\mathcal{A}}}
\newcommand{\Lm}{\ensuremath{\mathcal{L}}}

\newcommand{\mot}{\mathcal{M}^o_+(S)}
\begin{document}
\title{Relative motives and the theory of pseudo-finite fields}
\author[J. Nicaise]{Johannes Nicaise}
\address{Universit\'e Lille 1\\
Laboratoire Painlev\'e, CNRS - UMR 8524\\ Cit\'e Scientifique\\59 655 Villeneuve d'Ascq C\'edex \\
France} \email{johannes.nicaise@math.univ-lille1.fr}

\begin{abstract}
We generalize the motivic incarnation morphism from the theory of
arithmetic integration to the relative case, where we work over a
base variety $S$ over a field $k$ of characteristic zero. We
develop a theory of constructible effective Chow motives over $S$,
and we show how to associate a motive to any $S$-variety. We give
a geometric proof of relative quantifier elimination for
pseudo-finite fields, and we construct a morphism from the
Grothendieck ring of the theory of pseudo-finite fields over $S$,
to the tensor product of $\Q$ with the Grothendieck ring of
constructible effective Chow motives. This morphism yields a
motivic realization of parameterized arithmetic integrals.
Finally, we define relative arc and jet spaces, and the three
relative motivic Poincar\'e series.
\end{abstract}
\maketitle
\section{Introduction}
Let $k$ be a field of characteristic zero. Consider the
Grothendieck ring $K_0(PFF_k)$ of the theory of pseudo-finite
fields over $k$ (Section \ref{grothpff}), and denote by
$K_{0}^{mot}(Var_k)$ the image of the Grothendieck ring of
$k$-varieties in the Grothendieck ring of Chow motives,
 under the well-known
morphism $\chi_c:K_0(Var_k)\rightarrow K_0(CMot_k)$, see
\cite{GS}, \cite{GuiNav}, or \cite{Bitt} for a short construction
using weak factorization. Denef and Loeser \cite{DL} constructed a
generalized
 motivic Euler characteristic
$$\chi_{(c)}:K_{0}(PFF_k)\rightarrow K_{0}^{mot}(Var_k)\otimes \Q$$
 For this construction, it is important to understand the
structure of $K_{0}(PFF_k)$. The theory of quantifier elimination
for pseudo-finite fields \cite{FriJar,FriSac} states that
quantifiers can be eliminated if one adds some relations to the
language of rings, which have a geometric interpretation in terms
of Galois covers (Section \ref{Galois}). This interpretation
yields a construction for $\chi_{(c)}$. The morphism $\chi_{(c)}$
provides a concrete motivic realization of the theory of
arithmetic motivic integration \cite{DL}.

The goal of the present article, is to generalize this
construction to the relative case, where we work over an arbitrary
$k$-variety $S$ instead of over the field $k$, in order to obtain
a motivic incarnation of arithmetic integrals with parameters,
part of a work in progress by Cluckers and Loeser
\cite{ClLo,ClLo1,ClLo2}.

In Sections 2-4, we study the theory of pseudo-finite fields
(Definition \ref{pff}) over the base variety $S$. We define the
relative Grothendieck ring $K_0(PFF_S)$ in Section \ref{grothpff},
and we give a purely geometric proof of quantifier elimination for
ring formulas over $S$ w.r.t. the theory of pseudo-finite fields,
in terms of Galois formulas.

We briefly illustrate the concept of quantifier elimination by
means of a classic example. Let $\Lm$ be the first order language
 consisting of two binary function
symbols, a binary relation symbol, and two constant symbols, to
denote $+$, $.\,$, $<$, $0$ and $1$, respectively. These symbols
are interpreted in a structure with universe $\R$ in the obvious
way. By abuse of notation, we denote this structure again by $\R$.

Tarski proved that $\R$ allows elimination of quantifiers, meaning
that
%
there exists, for each formula $\psi(v_{1},\ldots,v_{m})$ in
$\Lm$, a formula $\varphi(v_{1},\ldots,v_{m})$ in $\Lm$, without
quantifiers, such that the set of $m$-tuples $(x_{1},\ldots
,x_{m})\in \R^{m}$ for which $\psi(x_{1},\ldots,x_{m})$ holds in
$\R$, is exactly the same as the set of $m$-tuples making the
formula $\varphi$ true in $\R$. In other words: every definable
subset of $\R^{m}$ can be described by a quantifier-free formula.
This is a very important property: it helps us to understand the
structure of definable sets, since quantifier-free formulas are
(in general) more transparent than formulas with quantifiers. The
key is to find a language describing a particular structure that
is sufficiently rich to allow elimination of quantifiers, but
sufficiently simple to keep quantifier-free definable sets
transparent.

Fried and Jarden developed in \cite{FriJar} a quantifier
elimination procedure for the theory of pseudo-finite fields by
introducing so-called Galois formulas (defined in Section
\ref{Galois}). Theorem \ref{rel} gives a purely geometric proof of
this result in the relative setting. It states that any ring
formula over $S$ is equivalent to a quantifier-free Galois
formula, where equivalent means that they define the same sets
when interpreted over an $M$-valued point of $S$, with $M$ a
pseudo-finite field.

To define a relative analogue of the Grothendieck ring
$K_0^{mot}(Var_k)$, we need an appropriate theory of motives over
an arbitrary base variety $S$ over $k$. The category $CMot_S$ of
constructible effective Chow motives over $S$ is constructed in
Section \ref{conchow} as a direct limit of
$$\prod_{S_i\in \mathscr{S}}\mathcal{M}^o_+(S_i)$$
where $\mathscr{S}$ runs over the finite stratifications of $S$
into smooth, irreducible locally closed subsets, and
$\mathcal{M}^o_+(S_i)$ is the category of effective Chow motives
over $S_i$ defined in \cite{deninger-murre}.

For any finite group $G$, we construct a functor from the category
$(G,Var_{S,c})$ of $S$-varieties with good $G$-action (with proper
morphisms), to  the homotopy category
$Ho(G,\mathcal{C}^b(CMot_S))$ of the additive category
$(G,CMot_S)$ of constructible motives with $G$-action. The
construction is based on a generalization of the extension
principle for cohomological functors \cite{GuiNav}, using the
existence of constructible resolution of singularities over $S$
(Proposition \ref{resolution}) and the existence of a split exact
blow-up sequence for relative effective Chow motives (Theorem
\ref{split}). In particular, we obtain a morphism of Grothendieck
rings $$\chi_c:K_0(Var_S)\rightarrow K_0(CMot_S)$$ whose image we
denote by $K_0^{mot}(Var_S)$.

We define induction and restriction functors on
$Ho(G,\mathcal{C}^b(CMot_S))$ with respect to morphisms of finite
groups, and we prove their main properties. The most important
feature of this formalism is the fact that the motive of a
quotient variety $X/G$ coincides with the $G$-invariant part of
the motive of $X$ (Theorem \ref{quot}). Using the induction and
restriction functors, we define a character decomposition of the
motive of an object in $(G,Var_S)$. This decomposition satisfies
certain Frobenius reciprocity properties (Lemma \ref{welldef3} and
Lemma \ref{welldef2}).

Now we are ready to construct the motivic realization morphism
$$\chi_{(c)}:K_0(PFF_S)\rightarrow K_0^{mot}(Var_S)\otimes \Q$$
Quantifier elimination over pseudo-finite fields, and the
geometric nature of Galois formulas, suggest a natural definition
of $\chi_{(c)}$. The subtle part is the proof that equivalent
formulas (i.e. formulas for which there exists a definable
bijection between the sets they define over pseudo-finite fields)
define the same motive; this follows from the Frobenius
reciprocity properties of the character decomposition.

Finally, we introduce relative arc and jet spaces. If $X$ is a
separated scheme of finite type over $S$, the $n$-th jet scheme
$\mathcal{L}_n(X/S)$ is a separated $S$-scheme of finite type that
parametrizes $(S\times_k k[t]/(t^{n+1}))$-valued sections on $X$.
The arc scheme $\mathcal{L}(X/S)$ is defined as a projective limit
of the jet schemes, and is endowed with natural projection
morphisms $$\pi_n:\mathcal{L}(X/S)\rightarrow \mathcal{L}_n(X/S)$$
We define the relative motivic Igusa Poincar\'e series (measuring
the jet schemes $\mathcal{L}_n(X/S)$), as well as the relative
geometric Poincar\'e series (measuring the projections
$\pi_n\mathcal{L}(X/S)$) and the arithmetic Poincar\'e series
(measuring the projections $\pi_n(\mathcal{L}(X/S)(K))$ for
pseudo-finite fields $K$). To show that the latter two are
well-defined, we establish a uniform version of Greenberg's
Theorem \cite{Gr}, based on a short new proof of the absolute
case. This result guarantees that the projections in the
definition are measurable (i.e. constructible, resp. definable by
a ring formula).

\vspace{8pt}

Let us give an overview of the results in this paper. In Section
\ref{pffield}, we formulate some basic results concerning Galois
covers and pseudo-finite fields, and we introduce the relative
Grothendieck ring of the theory of pseudo-finite fields. In
Section \ref{Galois}, we define relative Galois stratifications,
as well as some basic operations on these stratifications. Section
\ref{relelimination} is devoted to a geometric proof of the
elimination theorem (Theorem \ref{rel}).

In Section \ref{constr}, we prove the existence of constructible
equivariant resolution of singularities for varieties over $S$
with an action of a finite group, and a constructible equivariant
Hironaka-Chow Lemma. Section \ref{motives} is devoted to the
theory of constructible motives over $S$. Sections \ref{eff} to
\ref{blowup} contain some important properties of relative
effective Chow motives over a smooth base, in particular the
existence of a split exact blow-up sequence (Theorem \ref{split}).
In Section \ref{conchow}, we define our category of constructible
motives over $S$. The generalized extension criterion for
cohomological functors is stated in Section \ref{extension}, and
applied in Section \ref{groupmot} to associate a constructible
motive with $G$-action to any $S$-variety $X$ with good
$G$-action.

We define a character decomposition for constructible motives in
Section \ref{group}, and we prove its main properties. A crucial
result in this context is the fact that the motive of a quotient
$X/G$ coincides with the $G$-invariant part of the motive of $X$
(Theorem \ref{quot}). In Section \ref{relative}, Theorem
\ref{morph}, we construct a motivic incarnation morphism
$$\chi_{(c)}:K_0(PFF_S)\rightarrow K_0^{mot}(Var_S)\otimes \Q$$
for the relative theory of pseudo-finite fields.

Finally, in Section \ref{series}, we introduce relative arc and
jet spaces, and
 we prove a uniform version of Greenberg's Theorem,
which allows us to define the three relative motivic Poincar\'e
series.

\vspace{9pt}

 This paper borrows basic notions and
techniques from first order logic and model theory. A nice
introduction to these fields is given in \cite{marker}. In the
absolute case $S=\mathrm{Spec}\,k$, we recover the construction by
Denef and Loeser \cite{DL}. In our proofs, however, we avoided any
passage to finite fields, and hence any reference to
ultraproducts, merely using the geometric-arithmetic properties of
pseudo-finite fields established in Section \ref{pffield}. Our
proof of the quantifier elimination property in Theorem \ref{rel}
relies entirely on the existence of a certain short exact sequence
of algebraic fundamental groups. We hope this makes the
presentation more accessible to non-model theorists.

Throughout the paper, $k$ denotes a field of characteristic zero.
A variety over a scheme $S$ is a separated reduced scheme, of
finite type over $S$. For any scheme $S$, we denote by $S_{red}$
the underlying reduced scheme.

\section{Pseudo-finite fields}\label{pffield}
\subsection{Pseudo-finite fields and Galois covers}
\begin{definition}\label{pff}
A pseudo-algebraically closed field (PAC-field) is a field $M$,
such that every geometrically irreducible $M$-variety has an
$M$-rational point.

 A pseudo-finite field $M$ is an
infinite perfect PAC-field, which has exactly one field extension
of degree $n$, for every integer $n>0$, in a fixed algebraic
closure of $M$.
\end{definition}

Every field can be embedded in a pseudo-finite field
\cite[23.1.2]{FriJar}.
 An algebraic extension of a PAC field is PAC, by \cite[11.2.5]{FriJar}. As a consequence,
 a finite extension of a pseudo-finite field is pseudo-finite.
 Every finite extension of a pseudo-finite field is cyclic.
Ax \cite{Ax} proved that two ring formulas over $\Q$ are
equivalent when interpreted in $\mathbb{F}_p$, for all
sufficiently large primes $p$, if and only if they are equivalent
when interpreted in $K$, for all pseudo-finite fields $K$
containing $\Q$ \footnote{A ring formula over a
 field $k$ is a logical formula $\varphi$ built from Boolean combinations
 of polynomial equalities over $k$, and quantifiers.}. Hence, pseudo-finite fields
allow you to solve certain problems concerning finite fields in
characteristic zero.

\begin{definition}[Galois cover]
A Galois cover $h:Y\rightarrow X$ (also denoted by $Y/X$) is an
\'etale morphism of separated, integral, normal schemes,
satisfying the following property: there exists a finite group
$G$, acting faithfully on $Y$, and $h:Y\rightarrow X$ is a
quotient map for this action \cite[V.1]{sga1}.

The group $G$ (the group of $X$-automorphisms of $Y$) is called
the Galois group of the Galois cover $Y/X$, and is denoted by
$G(Y/X)$.
\end{definition}

\begin{definition}[Decomposition group]
Let $h:Y\rightarrow X$ be a Galois cover, and let $Z\rightarrow X$
be a morphism of separated, integral, normal schemes. Let $Z'$ be
a connected component of $Y\times_X Z$. The Galois group $G(Y/X)$
acts on $Y\times_X Z$. The decomposition group $D_{Y/X}(Z'/Z)$ of
$Z$ at $Z'$ with respect to $Y/X$, is the subgroup of elements of
$G(Y/X)$ which map the component $Z'$ to itself. It coincides with
the decomposition subgroup $D_{Y/X}(\eta_{Z'}/\eta_Z)$, where
$\eta_{Z'}$ and $\eta_Z$ are the generic points of $Z'$, resp.
$Z$.

The morphism $Z'\rightarrow Z$ is a Galois cover, and the
decomposition group of $Z$ at $Z'$ is canonically isomorphic to
the Galois group $G(Z'/Z)=G(\eta_{Z'}/\eta_Z)$. As $Z'$ runs
through the connected components of $Y\times_X Z$, the
decomposition group of $Z$ at $Z'$ runs through its conjugation
class in $G(Y/X)$. We call this conjugation class the
decomposition class
 of $Z$ (w.r.t. $Y/X$), and denote it by $C_{Y/X}(Z)$.
\end{definition}

 We list three easy properties of the
decomposition group for later use.

\begin{lemma}\label{decomp0}
Let $Y/ X$ be a Galois cover, let $M$ be a field, and let $x$ be
an $M$-valued point on $X$. Let $x'$ be a connected component of
$Y\times_X x$. Then $D_{Y/X}(x'/x)$ is the smallest subgroup $H$
of $G(Y/X)$, such that the point $x\rightarrow X$ lifts to a point
$x\rightarrow Y/H$, making the diagram $$\begin{CD} Y@<<< x'
\\@VVV @VVV
\\ Y/H@<<< x
\end{CD}$$
commute.
\end{lemma}
%

\begin{lemma}\label{decomp1}
Let $Y/X$ be Galois cover, let $V$ be a separated, normal,
integral scheme over $X$, and let $W$ be any connected component
of $V\times_X Y$. Let $A$ be a separated, normal, integral scheme
over $V$, and let $B$ be any connected component of $A\times_V W$.
If we denote by $B'$ the image of $B$ in $A\times_X Y$, then
 $D_{Y/X}(B'/A)=D_{W/V}(B/A)$
(where we view both sides as subgroups of $G(Y/X)$).
\end{lemma}
\begin{lemma}\label{decompb}
Let $Z/X$ be Galois cover, let $H$ be a normal subgroup of the
Galois group $G(Z/X)$, and let $Y$ be the quotient $Z/H$. Let $V$
be a separated, integral, normal scheme over $X$, let $V_Z$ be any
connected component of $Z\times_X V$, and let $V_Y$ be its image
in $Y\times_X V$. Then $D_{Y/X}(V_Y/V)$ is the image of
$D_{Z/X}(V_Z/V)$ under the restriction map $G(Z/X)\rightarrow
G(Y/X)$.
\end{lemma}

 A valuable property of pseudo-finite fields
 is given in the following lemma; see
\cite[24.1.4]{FriJar} (pseudo-finite fields are Frobenius fields).
\begin{lemma}\label{decomposition}
Let $Y\rightarrow X$ be a Galois cover, where $X$ and $Y$ are
varieties over a pseudo-finite field $M$ of characteristic zero,
and $X$ is
 geometrically irreducible over $M$. Denote by $N$ be the algebraic closure of $M$ in
the function field $k(Y)$. If $H$ is a cyclic subgroup of the
Galois group $G(Y/X)$, such that the image of $H$ under the
restriction morphism $G(Y/X)\rightarrow G(N/M)$ is the whole
Galois group $G(N/M)$, then there exist an $M$-rational point $x$
on $X$ such that $H$ belongs to $C_{Y/X}(x)$.
\end{lemma}

This lemma can be seen as an ``optimal lifting property'': for
each reasonable candidate, you find a rational point $x$ on $X$
with the desired decomposition group.

%

\begin{lemma}\label{decomp2}
Let $Y\rightarrow X$ a Galois cover, with $X$ and $Y$
 varieties over $k$.
If $\,C=<\!\sigma\!>$ is a cyclic subgroup of the Galois group
 $G(Y/X)$,
 then there exists
a closed $M$-valued point $x$ on $X$, with $M$ a pseudo-finite
field,
 such that $C$ belongs to $C_{Y/X}(x)$. Moreover, these points are
 dense in $X$.
\end{lemma}
\begin{proof}
If $U$ is a non-empty open subscheme of $X$, then $Y\times_X U/U$
is a new Galois cover with the same Galois group. Hence, if there
exist one such $x$, these points are dense in $X$.

Let $k'$ be the algebraic closure of $k$ in the function field
$k(X)$ of $X$. There exists a nonempty open subset $U$ of $X$,
such that $U$ is defined over $k'$. Hence, we may assume that
$k=k'$.

Let $K$ be the algebraic closure of $k$ in the function field
$k(Y)$, let $\sigma'$ be the image of $\sigma$ in $G(K/k)$, and
denote by $K^{\sigma'}$ its fixed field. Since $G(K/K^{\sigma'})$
is cyclic, there exists a pseudo-finite extension $M$ of
$K^{\sigma'}$, which is linearly disjoint from $K$ over
$K^{\sigma'}$, by \cite[23.1.1]{FriJar}.

The embedding of $K^{\sigma'}$ in $K$ defines a connected
component $Y'$ of $Y\times_k K^{\sigma'}$. The cover $Y'/X\times_k
K^{\sigma'}$ is Galois, and its Galois group is the decomposition
group of $X\times_k K^{\sigma'}$ at $Y'$ w.r.t. $Y/X$. It contains
the element $\sigma$. Since $M$ is linearly disjoint from $K$ over
$K^{\sigma'}$, the space $Y'\times_{K^{\sigma'}}M$ is connected,
and the Galois group $G(Y'\times_{K^{\sigma'}}M/X\times_k M)$ is
isomorphic to $G(Y'/X\times_k K^{\sigma'})$. In particular, it
contains $\sigma$ (if we identify
$G(Y'\times_{K^{\sigma'}}M/X\times_k M)$ with its image in
$G(Y/X)$).

By Lemma \ref{decomposition}, we can find a $M$-rational point $x$
on $X\times_k M$, whose decomposition class w.r.t.
$Y'\times_{K^{\sigma'}}M/X\times_k M$ contains $C$. By Lemma
\ref{decomp1}, its decomposition class w.r.t. $Y/X$ contains $C$,
as well.
\end{proof}

\begin{remark} Lemma \ref{decomp2} can be seen as a (weak) pseudo-finite
version of  Artin-Chebotarev's Density Theorem \cite{Serre2}:
assume that $X$, $Y$ are irreducible schemes of finite type over
$\Z$, of dimension $\geq 1$, and that $Y/G=X$ for some finite
group $G$ acting on $Y$,
 such that $G$ operates freely on $Y$, and
faithfully on the residue field of the generic point of $Y$.
 If $C$ is any subset of $G$, stable under conjugation,
then the set of closed points $x$ of $X$ for which $F_x\in C$, has
Dirichlet density equal to $|C|/|G|$. Here $F_x$ is the Frobenius
element in $G(k(y)/k(x))$, where $y$ is some closed point of $Y$,
lying over $x$.
\end{remark}

\subsection{The Grothendieck ring of the theory of pseudo-finite
fields}\label{grothpff} Let $S=\mathrm{Spec}\,R$ be an affine
scheme. We consider the relative Grothendieck ring $K_{0}(PFF_R)$
of the theory of pseudo-finite fields over $R$. As an abelian
group, it is generated by classes $[\varphi]$, where $\varphi$ is
a ring formula over $R$, which are subject to the relations
\begin{equation}\label{relation}
[\varphi_1\vee\varphi_2]=[\varphi_1]+[\varphi_2]-[\varphi_1 \wedge
\varphi_2], \end{equation}
 whenever $\varphi_1$ and $\varphi_2$
have the same free variables. For any pair of ring formulas
$\varphi_1$, $\varphi_2$ over $R$, with free variables, we impose
the additional relation $[\varphi_1]=[\varphi_2]$, whenever there
exists a ring formula $\psi$ over $R$ such that, for any
pseudo-finite field $M$ and any element $x$ of $S(M)$, the
interpretation of $\psi$ over $M$ defines the graph of a bijection
between the tuples over $M$ satisfying $\varphi_1$, and those
satisfying $\varphi_2$. Ring formulas over $R$ are interpreted
over $M$ in the obvious way, via the ring morphism $R\rightarrow
M$ corresponding to $x$. We denote this equivalence relation on
ring formulas over $R$ by $\varphi_1\equiv_{S} \varphi_2$, and we
say that $\psi$ defines an equivalence between $\varphi_1$ and
$\varphi_2$ over $S$. If $\varphi$ is a ring formula without free
variables, we make $\varphi$ equivalent to
$\varphi':=\varphi\wedge(v=0)$, where $v$ is a (free) variable,
and we impose $[\varphi]=[\varphi']$.

\vspace{10pt}
\begin{example}{Definable bijections.}
 Let $k$ be a field of characteristic zero, and put $R=k[z]$, and $S=\mathrm{Spec}\,R$.
 Let $f_1(x_1)$ and
$f_2(x_2)$ be polynomials in one variable over $k$, and consider
the ring formulas $\varphi_1=``f_1(x_1)=z"$ and
$\varphi_2=``f_2(x_2)=z"$ over $R$. Suppose that, for each point
$z_0$ of $S$, and each field extension $M$ of $k(z_0)$, the number
of elements $x_1$ that satisfy $\varphi_1$ over $M$, equals the
number of elements $x_2$ that satisfy $\varphi_2$ over $M$.

A trivial example is the case where $f_1$ and $f_2$ are linearly
related, i.e. there exist $a,b$ in $R$, with $a$ a unit, such that
$f_2(x)=f_1(ax+b)$. The ring formula $\psi=``x_1=ax_2+b"$ defines
an equivalence between $\varphi_1$ and $\varphi_2$ over $S$.

On the other hand, suppose $deg(f_1)=deg(f_2)=n$. For each integer
$m\geq 0$, and for $i=1,2$, there exists a ring formula
$\eta_{i,m}$ over $R$, without free variables, such that, for each
point $z_0$ of $S=\mathrm{Spec}\,R$, and each field extension $M$
of $k(z_0)$, the formula $\eta_{i,m}$ is true over $M$, iff the
equation $f_i(x)=z_0$ has exactly $m$ distinct solutions over $M$.
Hence, we can build a ring formula $\eta$ over $R$, without free
variables, such that for each point $z_0$ of $S$, and any field
extension $M$ of $k(z_0)$, the formula $\eta$ is true over $M$,
iff the solutions of $f_1(x_1)=z_0$ over $M$ are in bijective
correspondence with the solutions of $f_2(x_2)=z_0$ over $M$.
However, this does not mean that we can define (the graph of) such
a bijection by means of a ring formula over $R$. So, even if
$\eta$ is true for all $z_0$ and all $M$, we cannot conclude that
$\varphi_1\equiv_S \varphi_2$.


\end{example}

\vspace{5pt} We emphasize that, in our definition of $\equiv_S$,
it is important to consider every point of $S$, not only the
closed points. For example, let $p$ be a prime number, and let $R$
be the local ring of $\mathrm{Spec}\,\Z$ at $(p)$, with closed
point $x$, and generic point $\eta$.
 Then
$(p.1=0)\equiv_{x}(1=1)$, while $(p.1=0)\not\equiv_{\eta}(1=1)$.

 However, when $R$
is a ring of finite type over $\Z$, and $\varphi_1, \varphi_2,
\psi$ are ring formulas over $R$,
 then $\psi$ defines an equivalence of  $\varphi_1$ and $\varphi_2$ over $S$,
 iff, for each closed point $x$
of $S$, the formula $\psi$ defines an equivalence of $\varphi_1$
and $\varphi_2$ over $x$.
 This follows from
Ax' characterization of pseudo-finite fields as ultraproducts of
finite fields \cite{Ax}.

In fact, we will show in Section \ref{relelimination}, Corollary
\ref{generic}, that the same holds when $R$ is a
 ring of
finite type over any field $k$ of characteristic zero, using
quantifier elimination.

\vspace{8pt}
 Ring multiplication in the Grothendieck
ring is induced by taking the conjunction of formulas in disjoint
sets of free variables: if $\varphi_1$ and $\varphi_2$ are ring
formulas over $R$, with disjoint sets of free variables, then we
put $[\varphi_1].[\varphi_2]:=[\varphi_1\wedge \varphi_2]$. This
operation is well-defined, and extends bilinearly to a ring
product on $K_0(PFF_R)$.

 In order to understand the ring
$K_0(PFF_R)$, we need a relative quantifier elimination procedure,
which is described in Section \ref{relelimination}.

If $S$ is, more generally, a Noetherian scheme,
 we define $K_0(PFF_S)$ as the ring of global sections
of the unique Zariski-sheaf $\mathcal{F}(S)$ of rings on $S$, such
that, for each affine open subscheme $U=\mathrm{Spec}\,V$ of $S$,
the ring of sections $\mathcal{F}(S)(U)$ is equal to $K_0(PFF_V)$.
If $S=\coprod_{i}U_i$ is a finite stratification of $S$ into
locally closed affine subschemes $U_i$, the sheaf $\mathcal{F}(S)$
equals $\prod_i j^{(i)}_{!}\mathcal{F}(U_i)$, where $j^{(i)}$ is
the inclusion of $U_i$ in $S$. This is a consequence of relation
(\ref{relation}), since, for each $f$ in $R$, and each ring
formula $\varphi$ over $R$, we have $[\varphi]=[\varphi\wedge
f=0]+[\varphi\wedge f\neq 0]$.

Similarly, a ring formula $\varphi$ of $S$ consists of a finite
stratification $S=\coprod_{i}U_i$ of $S$ into locally closed
affine subschemes $U_i=\mathrm{Spec}\,R_i$, and a ring formula
$\varphi_i$ over $R_i$ for each $i$. A ring formula $\varphi$ over
$S$ defines a class $[\varphi]$ in $K_0(PFF_S)$. A quantifier-free
ring formula $\varphi$ with $m$ free variables defines a
constructible subset of $\A^{m}_S$. For any ring formula $\varphi$
over $S$, in $m$ free variables, for any point $x$ on $S$, and for
any pseudo-finite field extension $M$ of $k(x)$, we will denote by
$Z(\varphi,x,M)$, the subset of $M^{m}$ consisting of the tuples
satisfying the interpretation of $\varphi$ over $M$.

If $T$ is another Noetherian scheme, endowed with a morphism
$f:T\rightarrow S$,
 it is clear how to pull back a ring formula
over $S$ to a ring formula over $T$.

\begin{lemma}\label{pullback}
A morphism $f:T\rightarrow S$ induces a pull-back morphism
$$f^{*}:K_0(PFF_S)\rightarrow K_0(PFF_T)$$
\end{lemma}
\begin{proof}
If $\psi,\varphi_1, \varphi_2$ are ring formulas over $S$, such
that $\psi$ defines
 an equivalence between $\varphi_1$
and $\varphi_2$, the pull-back of $\psi$ defines an equivalence
between the pull-backs of $\varphi_1$ and $\varphi_2$.
\end{proof}

%

\section{Galois stratifications}\label{Galois}
\subsection{Definitions} We start by recalling the concepts of Galois
stratification and Galois formula; see \cite[\S 30]{FriJar} and
\cite[\S 2]{DL}.

 We fix a field $k$ of characteristic zero, as well as an irreducible
 $k$-variety $S$,
 and we will work in the
 category $(Var_S)$ of varieties over $S$.

\begin{definition}
 Let $X_S$ be a
variety over $S$.
  A normal stratification
$$<X_S,C_i/A_i\,|\,i\in I>\,$$
of $X_S$, is a partition of $X_S$ into a finite set of integral
and normal locally closed $S$-subvarieties $A_i$, each equipped
with a Galois cover $C_i\rightarrow A_i$.

A Galois stratification
$$\mathcal{A}=<X_S,C_i/A_i,Con(A_i)\,|\,i\in I>\,$$
is a normal stratification $<X_S,C_i/A_i\,|\,i\in I>$ with, for
each $i\in I$, a family $Con(A_i)$ of cyclic subgroups of the
Galois group $G(C_i/A_i)$, such that $Con(A_i)$ is stable under
conjugation. We call $Con(A_i)$ a conjugation domain for the cover
$C_i/A_i$. The support of $\mathcal{A}$ is the union of strata
with non-empty conjugation domain.
\end{definition}


For each point $x$ of the base scheme $S$, and for each
$S$-variety $Z_S$, we denote by $Z_x$ the fiber of $Z_S$ over $x$,
endowed with its reduced structure.
Let $\mathcal{A}$ be
a Galois stratification of $X_S$, let $x$ be a point of $S$, let
$M$ be a field extension of $k(x)$, and let $a$ be an $M$-valued
point of $X_x$, belonging to a stratum $A_{i,x}$. We put
$\Ar(C_i/A_i,x,a):=C_{C_i/A_i}(a)$. We write $\Ar(a)\subset
Con(\mathcal{A})$ for $\Ar(C_i/A_i,x,a)\subset Con(A_i)$. To these
data, we associate a subset $Z$ of $X_S(M)$ as follows:
$$Z(\mathcal{A},x,M)=\{a\in X_{x}(M)\,|\,\Ar(a)\subset
Con(\mathcal{A})\}\,$$ By Lemma \ref{decomp2}, the set
$Z(\mathcal{A},x,M)$ is empty for all pseudo-finite field
extensions $M$ of $k(x)$, iff the fiber of the support of
$\mathcal{A}$ over $x$ is empty.

Let $\mathcal{A}=<\A_{S}^{m+n},C_i/A_i,Con(A_i)\,|i\in I>$ be a
Galois stratification of $\A^{m+n}_{S}$, and let $Q_1,\ldots,Q_m$
be quantifiers. We denote by $\theta$, or by $\theta(\mathbf{Y})$,
the formal expression
$$(Q_1X_1)\ldots(Q_mX_m)[\Ar(\mathbf{X},\mathbf{Y})\subset
Con(\mathcal{A})]\,,$$ where $\mathbf{X}=(X_1,\ldots,X_m)$, and
$\mathbf{Y}=(Y_1,\ldots,Y_n)$. We call $\theta(\mathbf{Y})$ a
Galois formula over $S$ in the free variables $\mathbf{Y}$.

To a Galois formula $\theta$, to a point $x$ of $S$, and to a
field extension $M$ of $k(x)$, we associate the set
$$Z(\theta,x,M)=\{\mathbf{b}=(b_1,\ldots,b_n)\in M^{n}\,|\,(Q_1a_1)\ldots(Q_ma_m)\Ar(\mathbf{a},\mathbf{b})\subset
Con(\mathcal{A})\}$$ where the quantifiers $Q_1a_1,\ldots,Q_ma_m$
run over $M$.

\subsection{Ring formulas as Galois formulas}
It is easy to see how ring formulas over $S$ can be rewritten as
Galois formulas. Working locally, we may suppose
$S=\mathrm{Spec}\,R$. Let $\varphi(\mathbf{Y})$ be a formula in
the first order language of rings with coefficients in $R$, and in
the free variables $\mathbf{Y}$.
 Writing $\varphi$ in prenex normal form, we obtain a
formula
$$(Q_1X_1)\ldots(Q_mX_m)[\bigvee_{i=1}^{k}\bigwedge_{j=1}^{l}
f_{i,j}(\mathbf{X},\mathbf{Y})=0\wedge
g_{i,j}(\mathbf{X},\mathbf{Y})\neq 0]\,$$ with $f_{i,j}$ and
$g_{i,j}$ in $R[\mathbf{X},\mathbf{Y}]$. The formula between
brackets defines a constructible subset $W$ of $\A^{m+n}_{S}$, and
we can always find a stratification of $\A_{S}^{m+n}$ into
finitely many locally closed, integral, normal $S$-subvarieties
$A_i$, such that each stratum $A_i$ is either contained in $W$, or
in its complement. For each $i$, we take the Galois cover
$C_i\rightarrow A_i$ to be the identity, and we define $Con(A_i)$
as the family containing only the trivial group if $A_i$ is
contained in $W$, and as the empty family otherwise. In this way,
we obtain a Galois formula $\theta$, satisfying
$Z(\theta,x,M)=Z(\varphi,x,M)$ for each point $x$ of $S$, and each
field extension $M$ of $k(x)$.

\subsection{Galois formulas as ring formulas}\label{galoisring}
Now, we show how Galois formula $\theta$ can be rewritten as ring
formula $\varphi$, such that $Z(\theta,x,M)=Z(\varphi,x,M)$ for
each point $x$ of $S$, and each field extension $M$ of $k(x)$. We
may suppose that $S=\mathrm{Spec}\,R$. It is sufficient to
construct, for any integral, normal, locally closed subset $A$ of
$\A^m_S$, for any Galois cover $C/A$, and for any cyclic subgroup
$H$ of $G(C/A)$, a ring formula $\varphi_{C/A,H}$ with the
following property: for any point $x$ of $S$, and any field
extension $M$ of $k(x)$, the set $Z(\varphi_{C/A,H},x,M)$ is the
set of points $a$ in $A_x(M)$, such that $H\in C_{C/A}(a)$.

The locally closed subset $A$ can be defined by means of a ring
formula $\varphi_A$ over $R$ in $m$ free variables. Any \'etale
cover $D$ of $A$ can be defined by means of a ring formula
$\varphi_D$ over $R$ in $m+i$ free variables, for some integer
$i\geq 0$, such that the morphism $C\rightarrow A$ corresponds to
projection on the first $m$ coordinates. Now we can use Lemma
\ref{decomp0} to construct the formula $\varphi_{C/A,H}$.

Hence, using Galois formulas instead of ring formulas does not
alter the class of definable sets. The advantage of the new
formalism is the property of quantifier elimination, as we will
see in Theorem \ref{rel}.
\subsection{Refinement of Galois
stratifications} Let $X_S$ be an $S$-variety, and let
$\mathcal{A}=<X_S,C_i/A_i,Con(A_i)>$ be a Galois stratification.

Suppose that each $A_i$ is stratified into finitely many integral,
normal, locally closed subvarieties $A_{i,j}$. We will extend this
stratification to a Galois stratification
$\mathcal{A'}=<X_S,C_{i,j}/A_{i,j},Con(A_{i,j})>$, such that for
each point $x$ of $S$, and each field extension $M$ of $k(x)$,
$$Z(\Aa,x,M)=Z(\Aa',x,M)$$

Fix a stratum $A_{i,j}$, and let $C_{i,j}$ be a connected
component of $C_i\times_{A_i} A_{i,j}$. The projection
$C_{i,j}\rightarrow A_{i,j}$ is a Galois cover, and its Galois
group $G(C_{i,j}/A_{i,j})$ is the decomposition group of $C_i/A_i$
at $C_{i,j}$. We choose $Con(A_{i,j})$ to be the conjugation
domain consisting of the members of $Con(A_i)$ which are contained
in $G(C_{i,j}/A_{i,j})$. It follows from Lemma \ref{decomp1}, that
the resulting Galois cover $\mathcal{A}'$ satisfies our
requirements. We say that $\Aa'$ is induced from $\Aa$ by the
refinement $\{A_{i,j}\}$.

\subsection{Pulling back Galois
stratifications}\label{galoispback} The preceding refinement
procedure is a special case of a more general construction. Let
$f:Y\rightarrow X$ be any morphism of irreducible $S$-varieties,
 and let $\Aa$ be
a Galois stratification of $X$. We pull back $\mathcal{A}$ to a
Galois stratification $\mathcal{B}$ on $Y$ as follows: as
underlying stratification, we choose any stratification of $Y$
into finitely many integral, normal, locally closed subsets $B_j$,
 which is finer than the inverse image of the stratification of $X$.

Suppose $f$ maps $B_j$ into $A_i$, and let $D_j$ be any connected
component of $B_j\times_{A_i} C_i$. The projection $D_j\rightarrow
B_j$ is a Galois cover $D_j$ of $B_j$, and the Galois group
$G(D_j/B_j)$ is $D_{C_i/A_i}(D_j/B_j)$. We define the conjugation
domain $Con(B_j)$ as the set of members of $Con(A_i)$ which are
contained in $G(D_j/B_j)$.

We say that $\mathcal{B}$ is induced from $\mathcal{A}$ by $f$.
This construction is not canonical, but all choices of a
stratification correspond to equivalent Galois formulas, i.e. the
sets $Z(\mathcal{B},x,M)$ do not depend on any choices: for any
point $x$ on $S$, any field extension $M$ of $k(x)$, and any point
$a$ of $Y_x(M)$, we have
$$a\in Z(\mathcal{B},x,M)\mbox{ iff }f(a)\in Z(\Aa,x,M)$$ by Lemma \ref{decomp1}.

\subsection{Inflating a Galois stratification}
Let $C/A$ be a Galois cover, and let $Con(C/A)$ be a conjugation
domain for this cover. Suppose that $D/A$ is a Galois cover that
dominates $C/A$ (i.e. $D\rightarrow A$ factors through
$C\rightarrow A$). We define $Con(D/A)$ as the set of cyclic
subgroups of $G(D/A)$, whose image under the projection map
$G(D/A)\rightarrow G(C/A)$ belongs to $Con(C/A)$. This conjugation
domain has the following property: for any field $M$, and any
$M$-valued point $x$ on $A$, $C_{D/A}(x)\subset Con(D/A)$ iff
$C_{C/A}(x)\subset Con(C/A)$, by Lemma \ref{decompb}. We say that
$Con(D/A)$ is obtained by inflating the conjugation domain
$Con(C/A)$ to the cover $D/A$.

Let $X$ be an irreducible $S$-variety. We say that a Galois
stratification $$\Aa'=<X,C'_i/A'_i,Con(C'_i/A'_i)\,|\,i\in I>$$ is
obtained from a Galois stratification
$$\Aa=<X,C_i/A_i,Con(C_i/A_i)\,|\,i\in I>$$ by inflation, if
$A_i=A'_i$ for every $i\in I$ (modulo a permutation of the
indices), if $C'_i/A_i$ dominates $C_i/A_i$ for every $i\in I$,
and if the conjugation domain $Con(C'_i/A_i)$ is obtained from
$Con(C_i/A_i)$ by inflation. In this case,
$$Z(\Aa,x,M)=Z(\Aa',x,M)$$ for any point $x$ on $S$, and any field
extension $M$ of $k(x)$.
\subsection{Conjunction and disjunction of Galois stratifications}
Let $X$ be a variety over $S$. Consider Galois stratifications
$$\mathcal{A}=<\!Y,C_i/A_i,Con(A_i)\!>\mbox{ and }\mathcal{B}=<Y,D_j/B_j,Con(B_j)>$$
We will construct (non-canonical) Galois stratifications
$\mathcal{A}\vee \mathcal{B}$ and $\mathcal{A}\wedge \mathcal{B}$
with the following property: for any point $x$ on $S$, and any
field extension $M$ of $k(x)$,
\begin{eqnarray*}
Z(\mathcal{A}\vee \mathcal{B},x,M)&=&Z(\mathcal{A},x,M) \cup
Z(\mathcal{B},x,M) \\Z(\mathcal{A}\wedge
\mathcal{B},x,M)&=&Z(\mathcal{A},x,M) \cap Z(\mathcal{B},x,M)
\end{eqnarray*}
After a refinement, we may assume that $\mathcal{A}$ and
$\mathcal{B}$ have the same underlying stratification of $X$.
After an inflation process, we may even assume that the underlying
normal stratifications coincide. For any $i\in I$, we define the
conjugation domain for $C_i/A_i$ in the Galois stratification
$\mathcal{A}\vee \mathcal{B}$ (resp. $\mathcal{A}\wedge
\mathcal{B}$) as the union (resp. the intersection) of the
corresponding conjugation domains in $\mathcal{A}$ and
$\mathcal{B}$.
\subsection{Product of Galois stratifications}\label{galoisprod}
Let $X$ and $Y$ be varieties over $S$. Consider Galois
stratifications
$$\mathcal{A}=<\!X,C_i/A_i,Con(A_i)\!>\mbox{ and }\mathcal{B}=<Y,D_j/B_j,Con(B_j)>$$
We will construct a (non-canonical) Galois stratification
$\mathcal{A}\times_S \mathcal{B}$  of $X\times_S Y$, such that,
for any point $x$ on $S$, and any field extension $M$ of $k(x)$,
we have
$$Z(\mathcal{A}\times_S \mathcal{B},x,M)=Z(\Aa,x,M)\times
Z(\mathcal{B},x,M)\subset (X\times_S Y)_x(M)$$

Denote by $\pi_1$ and $\pi_2$ the projections of $X\times_S Y$ on
$X$, resp. $Y$. Pull back the Galois stratifications $\mathcal{A}$
and $\mathcal{B}$ to Galois stratifications $\mathcal{A}'$ and
$\mathcal{B}'$ of $X\times_S Y$, via $\pi_1$ and $\pi_2$. After a
refinement, we may assume that the underlying stratifications of
$\mathcal{A}'$ and $\mathcal{B}'$ coincide. Let $U$ be a stratum.
Let $C'/U$ and $D'/U$ be the Galois covers from $\mathcal{A}'$,
resp. $\mathcal{B}'$, with conjugation domains $Con(C'/U)$ and
$Con(D'/U)$. Let $V/U$ be any Galois cover dominating both $C'/U$
and $D'/U$. We define a conjugation domain $Con(V)$ as the
intersection of the inflations of $Con(C'/U)$ and $Con(D'/U)$ to
$V/U$. The constructed Galois cover has the desired property.

\subsection{Conjunction, disjunction and product of quantifier-free Galois formulas}
Fix integers $m,n\geq 0$, consider Galois stratifications
$$\mathcal{A}=<\!\A^m_S,C_i/A_i,Con(A_i)\!>\mbox{ and }\mathcal{B}=<\A^n_S,D_j/B_j,Con(B_j)>$$
and denote by $\theta_{\mathcal{A}}$ and $\theta_{\mathcal{B}}$
the corresponding Galois formulas.

If $m=n$, then we define
$\theta_{\mathcal{A}}\vee\theta_{\mathcal{B}}$ and
$\theta_{\mathcal{A}}\wedge\theta_{\mathcal{B}}$ as the Galois
formulas corresponding to the Galois stratifications
$\mathcal{A}\vee \mathcal{B}$, resp. $\mathcal{A}\wedge
\mathcal{B}$.

For any $m,n$, we define
$\theta_{\mathcal{A}}\times_S\theta_{\mathcal{B}}$ as the Galois
formula corresponding to the Galois stratifications
$\mathcal{A}\times_S \mathcal{B}$.

The Galois formulas are not canonically defined, but the sets
$Z(.,x,M)$ they define are independent of any choices, for any
point $x$ of $S$ and any field extension $M$ of $k(x)$.

\section{Quantifier elimination for Galois formulas}\label{relelimination}

\begin{theorem}[Relative quantifier elimination]\label{rel}
Let $S$ be a variety over a field $k$ of characteristic zero. Let
$\mathcal{A}$ be a Galois stratification of $\A^{m+n}_{S}$, and
let $\theta$ be a Galois formula
$$(Q_1X_1)\ldots(Q_mX_m)[\Ar(\mathbf{X},\mathbf{Y})\subset
Con(\mathcal{A})]\,$$ There exists a Galois stratification
$\mathcal{B}$ of $\A^{n}_{S}$ such that, for each point $x$ of
$S$, and each pseudo-finite field $M$ containing $k(x)$,
$$Z(\theta,x,M)=Z(\mathcal{B},x,M)\,$$
\end{theorem}

The absolute case $S=\mathrm{Spec}\,k$ was proven in
\cite{FriJar}. The relative case can immediately be reduced to the
absolute one, but we prefer to give a direct, purely geometric
proof, the novelty being the use of the exact sequence of
algebraic fundamental groups.

Theorem \ref{rel} follows from the elimination procedure in Lemma
\ref{existential} and Lemma \ref{universal} below. First, we need
an auxiliary result.

\begin{lemma}\label{connected}
Let $f:Y\rightarrow X$ be a proper, smooth morphism of irreducible
$k$-varieties, let $V$ be a relative strict normal crossing
divisor on $Y$ over $X$ \footnote{The definition of a relative
strict normal crossing divisor is recalled in Definition
\ref{relnc}.}, put $U=Y\setminus V$, and let $g:Z\rightarrow U$ be
a finite \'etale morphism of irreducible $k$-varieties. If $k(X)$
is algebraically closed in $k(Z)$, then the fibers of $h:=f\circ
g$ are geometrically irreducible.
\end{lemma}
\begin{proof}
Let $z$ be any geometric point on $U$, and put $x=f\circ z$. We
denote by $U_x$ and $Z_x$ the fibers of $U$, resp. $Z$, over $x$.

 By \cite{sga1}, Expos\'e XIII, Lemme 4.1 and
Exemples 4.4, there exists an exact homotopy sequence for the
morphism $f|_U:U\rightarrow X$ of the form
$$\begin{CD}\pi_1(U_x,z)@>f_1>>
\pi_1(U,z)@>f_2>> \pi_1(X,x)@>>> 1\end{CD}$$

 If $Z_x$ were not
connected, we could find a non-trivial \'etale cover $E$ of $U$,
dominated by $Z$, such that $E$ admits a section over $U_x$: it
suffices to take the quotient of $Z$ by the decomposition subgroup
of $U_x$ at any connected component of $Z_x$ w.r.t. $Z/U$.

By \cite[5.2.5]{Murre}, the fact that $ker(f_2)$ is contained in
$Im(f_1)$, has the following geometric interpretation: for any
connected \'etale cover $E$ of $U$, such that the restriction over
$U_x$ admits a section, there exists an \'etale cover $E'$ of $X$,
and a $U$-morphism from a connected component $C$ of $E'\times_X
U$ to
 $E$.

 However, since $k(X)$ is algebraically closed
in $k(Z)$, any such $E$ is trivial. Hence, $Z_x$ is connected.
\end{proof}

\begin{lemma}[Elimination of an existential quantifier]\label{existential}
 For every Galois stratification $\Aa$ of $\A^{m+1}_{S}$, there exists a Galois stratification $\mathcal{B}$ of $\A^{m}_{S}$, such that,
for each point $x$ of $S$, for each pseudo-finite field extension
$M$ of $k(x)$, and for each point $b$ of $\A^{m}_{x}(M)$, the
condition $\Ar(b)\subset Con(\mathcal{B})$ is equivalent to
$$(\exists a\in \A^{m+1}_{S}(M))(\pi(a)=b\wedge\Ar(a)\subset Con(\Aa))\,$$
where $\pi:\A^{m+1}_{S}\rightarrow \A^{m}_{S}$ is the projection
on the first $m$ coordinates.
\end{lemma}
\begin{proof}
Let $\Aa$ be a Galois stratification of $\A^{m+1}_{S}$.
 Refining
$$\mathcal{A}=<\A^{m+1}_S,C_i/A_i,Con(A_i)\,|\,i\in I>$$ if necessary, we may suppose
that there exists a finite stratification $\{B_j\}_{j\in J}$ of
$\A^m_S$ into normal, integral, locally closed subsets, such that
for any $i\in I$, there exists an index $j(i)\in J$ such that
$\pi:A_i\rightarrow \A^m_S$ factors through a smooth surjective
morphism $\pi:A_i\rightarrow B_{j(i)}$. This stratification is
constructed using generic smoothness and Noetherian induction.

For each $j\in J$, and each $i\in I$ with $j(i)=j$, we will
construct a Galois cover $D^i_j/B_j$ and a conjugation domain
$Con(D^i_j/B_j)$ with the following property: for any
pseudo-finite field $M$, and any $M$-valued point $b$ on $B_j$,
$C_{D_j^i/B_j}(b)\subset Con(D_j^i/B_j)$ iff there exists an
$M$-valued point $a$ on $A_i$ with $\pi(a)=b$ and
$C_{C_i/A_i}(a)\subset Con(A_i)$. Dominating the covers
$D^i_j/B_j$ by a common Galois cover $D_j/B_j$ and inflating the
conjugation domains, we may suppose that all the covers
$D^i_j/B_j$ coincide. Finally, we put $Con(B_j):=\cup_{i\in
I,j(i)=j}Con(D^i_j/B_j)$, and
$$\mathcal{B}:=<\A^{m}_S,D_j/B_j,Con(B_j)>$$
The Galois stratification $\mathcal{B}$ satisfies the
requirements.

We fix $j\in J$, and $i\in I$ with $j=j(i)$. We'll simply write
$C/A,\,D/B$ and $Con(B)$, instead of $C_i/A_i,\,D^i_j/B_j$ and
$Con(D^i_j/B_j)$.

 We will distinguish two cases: either the dimensions of
$A$ and $B$ agree, or dim$\,A=$\ dim$\,B+1$.

 \textit{Case 1:
dim\,$A$=dim\,$B$.} Stratifying $B$ if necessary, we may suppose
that $\pi:A\rightarrow B$ is \'etale and finite. Hence, $C$ is an
\'etale cover of $B$, and can be dominated by a Galois cover
$D\rightarrow B$. We define $Con(B)'$ as the inflation of $Con(A)$
to $D/A$, and $Con(B)$ as the smallest conjugation domain in
$G(D/A)$ containing $Con(B)'$. Here we consider $G(D/A)$ as a
subgroup of $G(D/B)$.

\textit{Case 2: dim\,$A$=dim\,$B+1$.} Then $\pi:A\rightarrow B$
factors through an open immersion $A\rightarrow \mathrm{Spec}\,
\mathcal{O}_{B}[t]$.

If we denote by $D$ the normalization of $B$ in the algebraic
closure of $k(B)$ in $k(C)$, then $D$ is a Galois cover of $B$.
Intuitively, this construction extracts the Galois action on the
base $B$ from the cover $C/A$. The cover $C/A$ factors through the
Galois cover $A\times_{B}D\rightarrow A$, and
\begin{equation}\label{iso} G(A\times_{B}D/ A)\cong
G(D/B)\end{equation} Let $Con(B)$ be the conjugation domain,
obtained by restricting the elements of the members of $Con(A)$ to
$G(D/B)$ via
$$G(C/A)\rightarrow G(A\times_{B}D/ A)\cong
G(D/B) $$

Let $x$ be a point of $S$, let $M$ be a field extension of $k(x)$,
and let $b$ be a point of $B_{x}(M)$. First, suppose that there
exists a point $a$ of $A(M)$, with $\pi(a)=b$, such that
$C_{C/A}(a)\subset Con(A)$. Since $C_{D/B}(b)$ is obtained from
$C_{C/A}(a)$ by restriction (by Lemma \ref{decompb} and the
isomorphism (\ref{iso})), $C_{D/B}(b)\subset Con(B)$.

Now suppose, conversely, that $C_{D/B}(b)\subset Con(B)$, and that
$M$ is pseudo-finite. Let $b'$ be a point on the
 fiber of $D/B$ over $b$. Its decomposition group is
an element $H_0$ of $Con(B)$. We denote its residue field $k(b')$
by $N$.


Put $F_{b'}=A\times_{B}b'$, and let $C_{b'}$ be the inverse image
 of $F_{b'}$ under the Galois cover $C/A\times_{B}D$. Stratifying the base $B$, we may
suppose that the complement of $A$ in $\mathrm{Spec}\,
\mathcal{O}_{B}[t]$ is either empty, or \'etale and finite over
$B$. This means, in particular, that $A\times_{B}D$ is the
complement in $D\times \Pro^{1}$ of a divisor with normal
crossings, relative to $D$.
By Lemma \ref{connected}, this implies that $C_{b'}$ is
geometrically connected over $k(b')$.  Hence, the decomposition
group $D_{C/A\times_B D}(C_{b'}/F_{b'})$ is the whole Galois group
$G(C/A\times_{B}D)$, i.e. the Galois action of $C/A\times_{B}D$ is
concentrated on the fibers.

By definition of $Con(B)$, there exists an element $H_1$ of
$Con(A)$, such that $H_0$ is the image of $H_1$ under the
restriction morphism $G(C/A)\rightarrow G(D/B)$. Let $F_b$ be the
fiber of $A$ over $b$, and consider the induced Galois cover
$C_{b'}/F_b$. Since $C_{b'}$ is geometrically connected over
$k(b')$, the algebraic closure of $M$ in $k(C_{b'})$ is $N$.
By Lemma \ref{decomposition}, there exists a point $a$ of $F_b(M)$
with decomposition group $H_1\in Con(A)$.
\end{proof}

\begin{lemma}[Elimination of a universal quantifier]\label{universal}
 For every Galois stratification $\Aa$ of $\A^{m+1}_{S}$, there exists a Galois stratification $\mathcal{B}$ of $\A^{m}_{S}$, such that,
for each point $x$ of $S$, for each pseudo-finite field extension
$M$ of $k(x)$, and for each point $b$ of $\A^{m}_{x}(M)$, the
condition $\Ar(b)\subset Con(\mathcal{B})$ is equivalent to
$$(\forall a\in \A^{m+1}_{S}(M))(\pi(a)=b\Rightarrow\Ar(a)\subset Con(\Aa))$$
\end{lemma}
\begin{proof}
We will deduce Lemma \ref{universal} from Lemma \ref{existential}.
Let $\Aa^{c}$ be the complementary Galois stratification of $\Aa$;
this means that $\Aa^{c}$ has the same underlying normal
stratification, but for each stratum $A_i$, the conjugation domain
$Con^{c}(A_i)$ associated to $A_i$ by $\Aa^{c}$ consists of the
cyclic subgroups of $G(C_i/A_i)$ that do not belong to $Con(A_i)$.
Lemma \ref{existential} produces a Galois stratification
$\mathcal{B}^{c}$ of $\A^{m}_{S}$; its complement $\mathcal{B}$
satisfies the requirements.
\end{proof}

\noindent Repeatedly applying Lemma \ref{existential} and Lemma
\ref{universal} proves Theorem \ref{rel}.

\begin{cor}\label{generic}
Let $S$ be an irreducible $k$-variety, and let
$\varphi_1,\varphi_2,\psi$ be ring formulas over $S$. Then $\psi$
defines an equivalence of $\varphi_1$ and $\varphi_2$ over the
generic point $\eta$ of $S$, iff $\psi$ defines an equivalence of
$\varphi_1$ and $\varphi_2$ over each closed point $x$ in a dense
open subset $U$ of $S$.
\end{cor}
\begin{proof}
The property '$\psi$ defines a bijection between $\varphi_1$ and
$\varphi_2$' is definable by a sentence $\theta$, i.e. a ring
formula over $S$ without free variables. Hence, it suffices to
show that $\theta$ is true in every pseudo-finite extension of
$k(\eta)$, iff it is true in every pseudo-finite extension of
$k(x)$, for each closed point $x$ in some dense open subset $U$ of
$S$.

 By the elimination theorem, we may suppose that
$\theta$ is a Galois formula without quantifiers and without free
variables. Let $U$ be the stratum of the corresponding Galois
stratification, containing $\eta$. Let $C/U$ be the associated
Galois cover, and $Con(U)$ its conjugation domain.

The interpretation of $\theta$ is true, in every pseudo-finite
extension of $k(\eta)$, iff $Con(U)$ contains the decomposition
classes of $C/U$ at all pseudo-finite extensions of $k(\eta)$. By
Lemma \ref{decomp2}, applied to the Galois cover $k(C)/k(\eta)$,
this is the case iff $Con(U)$ contains all cyclic subgroups of
$G(C/U)$. Again by Lemma \ref{decomp2}, this is, at its turn,
equivalent to saying that the interpretation of $\theta$ is true,
in every pseudo-finite extension of every closed point $x$ of $U$.
\end{proof}

\begin{cor}\label{generators}
As an Abelian group, the Grothendieck ring $K_0(PFF_S)$ is
generated by classes of the form $[\varphi_{Y/X,C}]$, where $Y/X$
is a Galois cover of affine normal irreducible $S$-varieties, and
$C$ is a cyclic subgroup of $G(Y/X)$ (see Section \ref{galoisring}
for the definition of the formula $\varphi_{Y/X,C}$).
\end{cor}
\begin{proof}
By the elimination theorem, $K_0(PFF_S)$ is generated by the
classes of (ring formulas corresponding to) quantifier-free Galois
formulas $\theta$. Writing the conjugation domain of $\theta$ as a
union of conjugation classes of cyclic groups, we see that the
class of $\theta$ in $K_0(PFF_S)$ can be written as a sum of
classes of the form $[\varphi_{Y/X,C}]$.
\end{proof}
\section{Constructible resolution of singularities}\label{constr}
\subsection{Varieties with good group action}
Let $S$ be a variety over a field $k$ of characteristic zero, and
let $G$ be a finite group. If $X$ is any variety over $S$, a good
$G$-action on $X$ is an action of $G$ on $X$, such that the
structural morphism $X\rightarrow S$ is equivariant (where $S$
carries the trivial action), and such that every orbit is
contained in an affine open subscheme of $X$ (this is automatic if
$X$ is quasi-projective over $S$). We define the Grothendieck
group $K_0^{G}(Var_S)$
as follows: start with the free abelian group on isomorphism
classes $[X]$ of $S$-varieties $X$ with good $G$-action, and
consider relations of the form $[X]=[Z]+[X\setminus Z]$, where $Z$
is a $G$-invariant closed subvariety with the induced $G$-action.
The fiber product over $S$ induces a ring product on
$K_0^{G}(Var_S)$. If $X$ is an $S$-variety, and $S_i$ is a locally
closed subset of $S$, we denote by $X_{S_i}$ the $S_i$-variety
$(X\times_S S_i)_{red}$.

For any morphism of $k$-varieties $f:S'\rightarrow S$, there is a
base-change morphism of rings $$f^*:K_0^G(Var_S)\rightarrow
K_0^G(Var_{S'})$$ For any morphism of finite groups $G\rightarrow
G'$, there is a forgetful morphism
$$K_0^{G'}(Var_S)\rightarrow
K_0^G(Var_{S})$$ If $G$ is the trivial group $\{e\}$, we write
$K_0(Var_S)$ instead of $K_0^G(Var_S)$.
\subsection{Resolution of singularities}

Let $S$ be a variety over a field $k$ of characteristic zero.

\begin{definition}
An admissible stratification is a finite stratification
$\mathscr{S}=\{S_1,\ldots,S_p\}$ of $S$ into smooth, irreducible,
quasi-projective, locally closed subvarieties $S_i$.
\end{definition}

\begin{definition}[\cite{sga1},XIII.2.1]\label{relnc}
Let $X$ be a smooth $S$-variety. A strict normal crossing divisor
relative to $S$ on $X$ is a divisor on $X$, with the following
properties:
\begin{itemize}
\item the prime components of $D$ are smooth over $S$, \item at
any point $x$ on $X$, these components are locally defined by
$f_i=0$, $i=1,\ldots,r$, with $f_i(x)=0$. The scheme
$V(f_1,\ldots,f_r)$ is smooth over $S$, of codimension $r$ in $X$.
\end{itemize}
\end{definition}

\begin{definition}
A resolution of singularities for a variety $X$ over $S$, is a
composition of blow-ups with $S$-smooth centers $h:X'\rightarrow
X$, such that $X'$ is smooth over $S$, and such that the
exceptional locus $E$ of $h$ is a strict normal crossing divisor
relative to $S$.
\end{definition}

If $S=\mathrm{Spec}\,k$, any variety $X$ over $S$ admits a
resolution of singularities, by Hironaka's famous result
\cite{hironaka}. This does not hold for arbitrary $S$. The
following result is a constructible version of resolution of
singularities over $S$.

\begin{prop}[Constructible resolution of singularities]\label{resolution}
For any variety $X$ over $S$, there exists an admissible
stratification $\mathscr{S}=\{S_i\,|\,i\in I\}$ of $S$, such that
for any $i\in I$, the variety $(X\times_S S_i)_{red}$ admits a
resolution of singularities over $S_i$.

If $G$ is a finite group, and $X$ carries a good $G$-action, then
all the blow-ups in the resolution can be chosen to have
$G$-closed centers, so that the resolution of singularities is
$G$-equivariant. Moreover, if $V$ is a $G$-closed, closed
subvariety of $X$, we can find a resolution by blow-ups of smooth
$G$-closed centers, such that the union of the exceptional locus
with the inverse image of $V$ is a strict normal crossing divisor
relative to $S$.
\end{prop}
\begin{proof}
It suffices to prove that there exists a non-empty open subvariety
$U$ of $S$, such that $X\times_S U$ admits a resolution over $U$
by blow-ups with $G$-closed centers. The theorem then follows by
 Noetherian induction on $S$.

By $G$-equivariant uniformization of ideals \cite[\S 7]{vmayor2},
the result holds for the fiber $X_\eta$ of $X$ over the generic
point $\eta$ of $S$. Hence, it suffices to prove the following
claims.

\textit{1. If $X_\eta$ is smooth over $\eta$, then there exists a
non-empty open subscheme $U$ of $X$, such that $X_U$ is smooth
over $U$.}

\noindent We may assume that $X$ and $S$ are affine. If the
Jacobian criterion for smoothness holds over $\eta$, it holds over
an open neighbourhood $U$ of $\eta$ in $S$.

\textit{2. If $Z_\eta$ is a closed smooth subvariety of $X_\eta$,
then there exists a non-empty open subscheme $U$ of $S$, and a
smooth closed subvariety $Z_U$ of $X_U:=X\times_S U$, such that
$Z_\eta$ is the fiber of $Z_U$ over $\eta$. The blow-up of
$X_\eta$ at $Z_\eta$ is the fiber over $\eta$ of the blow-up of
$X_U$ at $Z_U$.}

\noindent Working locally, we may assume that $S$ and $X$ are
affine. The equations defining $Z_\eta$ in $X_\eta$ extend to
regular functions on $X_U$, for some non-empty open subscheme $U$
of $X$, and these define a closed subvariety $Z_U$ with
$Z_\eta=Z_U\times_U \eta$. By point 1, after shrinking $U$, we may
assume that $Z_U$ is smooth over $U$. The last part of the claim
follows from flat base change for blow-ups \cite[8.1.12.c]{Liu}.

\textit{3. If $D$ is a divisor on $X$, and $D_\eta:=D\times_X
\eta$ is a strict normal crossing divisor on $X_\eta$, then there
exists a non-empty open subscheme $U$ of $S$ such that
$D_U:=D\times_X U$ is a strict normal crossing divisor relative to
$U$ on $X_U$.}

\noindent Apply point 1 to the components of $D$ and their
intersections.
\end{proof}

\begin{prop}[Constructible Chow-Hironaka Lemma]\label{chow}
Let $G$ be a finite group. Let $X,X'$ be smooth irreducible
$S$-varieties of relative dimension $d$, carrying a good
$G$-action, and let $Y$ be a closed subvariety of $X'$, such that
the dimensions of the fibers of $Y\rightarrow S$ are strictly
smaller than $d$. Let $h:X'\rightarrow X$ be a $G$-equivariant
proper birational morphism over $S$, such that $h$ is an
isomorphism over $X\setminus Y$.

There exists an admissible stratification
$\mathscr{S}=\{S_i\,|\,i\in I\}$ such that, for any $i\in I$, the
morphism $h_{S_i}:X'_{S_i}\rightarrow X_{S_i}$ can be dominated by
a composition of blow-ups with $S_i$-smooth, $G$-closed centers.
\end{prop}
\begin{proof}
By the equivariant Chow-Hironaka Lemma, the proposition holds if
$S$ is the spectrum of a field (see \cite{DL5}, Lemma A.2). Now we
can proceed as in the proof of Proposition \ref{resolution}.
\end{proof}

\begin{prop}[Constructible compactification]\label{compact}
Let $G$ be a finite group. Let $X$ be a smooth $S$-variety,
carrying a good $G$-action, and let $V$ be a $G$-closed, closed
subvariety. There exists an admissible stratification
$\mathscr{S}=\{S_i\,|\,i\in I\}$ such that, for any $i\in I$, we
can find a $G$-equivariant compactification
$X_{S_i}\hookrightarrow X'_{S_i}$, with $X'_{S_i}$ a smooth,
proper $S_i$-variety with good $G$-action, and such that the union
of $X'_{S_i}\setminus X_{S_i}$ with the closure of $V_{S_i}$ is a
strict normal crossing divisor relative to $S_i$.
\end{prop}
\begin{proof}
Starting with a $G$-equivariant compactification of $X$ over $S$,
we can use Proposition \ref{resolution} to obtain the result.
\end{proof}
\section{Constructible Chow motives over a base variety}\label{motives}
For the categorical language in this section (additive and
pseudo-abelian categories, tensor structures, $\ldots$) we refer
to \cite{Lev}.

\subsection{Pseudo-abelian categories}
An additive category $\mathcal{A}$ is called pseudo-abelian, if
all projectors (i.e. all idempotent endomorphisms) split. This
means that, for any projector $p$ on an object $A$, we can find an
isomorphism $A\cong B\oplus C$ such that $p$ corresponds to
$(Id,0)$ on $B\oplus C$. In this case, the natural morphism
$C\rightarrow A$ is a kernel for $p$, and $B\rightarrow A$ is an
image. For any additive category $\mathcal{A}$, the pseudo-abelian
envelope is an additive full embedding $\mathcal{A}\rightarrow
\mathcal{A}_{\sharp}$ defined by the following universal property:
$\mathcal{A}_{\sharp}$ is pseudo-abelian, and any additive functor
from $\mathcal{A}$ to a pseudo-abelian category $\mathcal{C}$
factors through an essentially unique additive functor from
$\mathcal{A}_{\sharp}$ to $\mathcal{C}$.

The pseudo-abelian envelope $\mathcal{A}_{\sharp}$ of an additive
category $\mathcal{A}$ is constructed by artificially adding
images for projectors (we get kernels for free, since the kernel
of a projector $p$ coincides with the image of $Id-p$). The
objects of $\mathcal{A}_{\sharp}$ are pairs $(A,p)$, where $A$ is
an object of $\mathcal{A}$ and $p$ is a projector on $A$. The
morphisms are given by
$$Hom_{\mathcal{A}_{\sharp}}((A,p),(B,q)):=q\circ Hom_{\mathcal{A}}(A,B)\circ p $$

\subsection{The Grothendieck ring of an additive category}
\begin{definition}[Grothendieck ring]\label{groth}
For any additive category $\mathcal{A}$, the Grothendieck group
$K_0(\mathcal{A})$ is the abelian group generated by the monoid
$([\mathcal{A}],\oplus)$, where $[\mathcal{A}]$ denotes the
set\footnote{We tacitly assume that the isomorphism classes of
objects of $\Aa$ form a set.} of isomorphism classes of objects of
$\mathcal{A}$, and $\oplus$ is the direct sum. We denote the class
of an object $A$ of $\mathcal{A}$ in $K_0(\mathcal{A})$ by $[A]$.

If $\mathcal{A}$ carries a tensor structure, the tensor product
$\otimes$ induces a multiplication on $K_0(\mathcal{A})$,
determined by $[A].[B]:=[A\otimes B]$. This multiplication turns
$K_0(\mathcal{A})$ into a ring: the Grothendieck ring associated
to $\mathcal{A}$.

An additive functor $F:\Aa\rightarrow \mathcal{B}$  respects
direct sums, and hence induces a morphism of Abelian groups
$K_0(F):K_0(\Aa)\rightarrow K_0(\mathcal{B})$. If $\Aa$ and
$\mathcal{B}$ carry tensor structures, and $F$ respects tensor
structures, $K_0(F)$ is a morphism of rings.
\end{definition}

\begin{remark}
This construction is not to be confused with the definition of the
Grothendieck group of an exact category (the free abelian group on
isomorphism classes, modulo the relation $[B]=[A]+[C]$ whenever
$0\rightarrow A\rightarrow B\rightarrow C\rightarrow 0$ is a short
exact sequence). For instance; consider the category
$\mathcal{A}b_{fg}$ of finitely generated abelian groups. Any such
group $G$ can be written canonically as a direct sum of $\Z^r$ for
some integer $r>0$, and a finite number of torsion groups
$\Z/(p^i)$, with $p$ a prime. The rank map
$$rk:\mathrm{Ob}(\Aa)\rightarrow \Z:G\mapsto r$$ induces an isomorphism from
the Grothendieck group of the abelian category $\Aa b_{fg}$ to
$\Z$: torsion parts are killed, by the existence of an exact
sequence $0\rightarrow \Z\rightarrow \Z\rightarrow
\Z/(n)\rightarrow 0$ for any integer $n>0$. However, if we
consider $\Aa b_{fg}$ merely as an additive category, we can write
any element $[G]$ of $K_0(\Aa b_{fg})$ uniquely as
$$r[\Z]+\sum_{p>0\,\mbox{prime}, i>0}n_{p,i}[\Z/(p^i)]$$ with
$r,n_{p,i}$ integers, and $n_{p,i}=0$ for almost all $p$ and $i$.
\end{remark}

\subsection{Effective Chow motives over a smooth base}\label{eff}
Let $k$ be a field of characteristic zero, and let $S$ be a smooth
irreducible quasi-projective variety over $k$. We briefly recall
the construction of the category $\mathcal{M}^{o}_{+}(S)$ of Chow
motives over $S$ (see \cite[\S 1]{deninger-murre}), and we prove
some elementary properties.

Denote by $\mathcal{V}_S$ the category of smooth, projective
varieties $X$ over $S$. If $X$ is irreducible, we denote the
relative dimension of $X$ over $S$ by $d(X/S)$. For any smooth
projective variety $X$ over $S$, and any integer $\alpha\geq 0$,
we denote by $CH^{\alpha}(X)$ the Chow group of algebraic cycles
of degree $\alpha$. It is constructed as follows: consider the
free abelian group on the set of irreducible subvarieties of $X$
of codimension $\alpha$, and take the quotient modulo rational
equivalence. We put
$CH^{\alpha}(X,\Q):=CH^{\alpha}(X)\otimes_{\Z}\Q$. For any pair of
integers $\alpha,\beta\geq 0$, there is a bilinear intersection
pairing $$CH^{\alpha}(X,\Q)\times CH^{\beta}(X,\Q)\rightarrow
CH^{\alpha+\beta}(X,\Q):(x,y)\mapsto x.y$$ which makes the graded
group $CH(X,\Q):=\oplus_{\alpha}CH^{\alpha}(X,\Q)$ into a ring.

If $X$, $Y$ and $Z$ are smooth, projective varieties over $S$,
with $X$ and $Y$ irreducible, we can construct a bilinear map
$$\circ\,:\,CH^{d(X/S)}(X\times_S Y,\Q)\times CH^{d(Y/S)}(Y\times_S Z,\Q)\rightarrow CH^{d(X/S)}(X\times_S Z,\Q)$$
as follows: a couple $(f,g)$ is mapped to
$$f\circ g:= (p_{13})_{*}(p_{12}^{*}(f).p_{23}^{*}(g))$$
where $p_{12},p_{13}, p_{23}$ are the projections of $X\times_S
Y\times_S Z$ on $X\times_S Y$, $X\times_S Z$, resp. $Y\times_S Z$,
and ``$.$'' is the intersection product in the Chow ring of
$X\times_S Y\times_S Z$.

Let $C\mathcal{V}_S^o$ be the category with the same objects as
$\mathcal{V}_S$, but with morphisms
$$Hom_{C\mathcal{V}_S^o}(X,Y):=\oplus_{X_i\in \mathcal{C}(X)}CH^{d(X_i/S)}(X_i\times_S Y,\Q)$$
where $\mathcal{C}(X)$ is the set of connected components of $X$.
The bilinear map $\circ$ constructed above defines the composition
law in $C\mathcal{V}_S^o$, turning it into an additive category
with $\Q$-linear structure; the direct sum $X\oplus Y$ is simply
the disjoint union $X\sqcup Y$. There is a canonical contravariant
functor
$$M:\mathcal{V}_S\rightarrow C\mathcal{V}^o_S$$ mapping a $S$-morphism
$f:Y\rightarrow X$ to the transpose of its graph.

The category $\mathcal{M}_{+}^o(S)$ of (effective) Chow motives
over $S$, is defined as the pseudo-abelian envelope of
$C\mathcal{V}_S^o$. The fiber product over $S$ induces a tensor
structure on $\mathcal{M}_{+}^o(S)$. The functor $M$ induces a
canonical contravariant functor
$$M:\mathcal{V}_S\rightarrow \mathcal{M}_{+}^o(S)$$
which we will denote by $M_S$ if we want to make the base
explicit. To simplify notation, we will sometimes denote the image
of a morphism of smooth projective $S$-varieties $f:X\rightarrow
Y$ by $f^*\in Hom_{\mot}(M(Y),M(X))$. This notation is not to be
confused with the base change functor defined below.

Applying Definition \ref{groth} to the additive tensor category
$\mathcal{M}_+^o(S)$, we obtain the Grothendieck ring of Chow
motives $K_0(\mathcal{M}_+^o(S))$.

%
%

If $f:S'\rightarrow S$ is a morphism of smooth, irreducible,
quasi-projective $k$-varieties, then the base change functor
$f^{*}:\mathcal{V}_S\rightarrow \mathcal{V}_{S'}$ induces a
$\Q$-linear base change functor
$$f^{*}:\mathcal{M}_{+}^o(S)\rightarrow \mathcal{M}_{+}^o(S')$$
compatible with the tensor structures, and a ring morphism
$$f^*:K_0(\mathcal{M}_+^o(S))\rightarrow K_0(\mathcal{M}_+^o(S'))$$

\begin{lemma}\label{forget}
If $f:S'\rightarrow S$ is a smooth, projective morphism of smooth,
irreducible, quasi-projective $k$-varieties, then the forgetful
functor $f_{*}:\mathcal{V}_{S'}\rightarrow \mathcal{V}_S$ induces
a $\Q$-linear forgetful functor
$$f_{*}:\mathcal{M}_{+}^o(S')\rightarrow \mathcal{M}_{+}^o(S)$$
(generally \textit{not} compatible with the tensor structures),
and a morphism of abelian groups
$$f_*:K_0(\mathcal{M}_+^o(S'))\rightarrow K_0(\mathcal{M}_+^o(S))$$
\end{lemma}
\begin{proof}
Let $X,Y$ be smooth and projective varieties over $S'$. The
canonical morphism $\pi:X\times_{S'}Y\rightarrow X\times_{S}Y$ is
proper, and induces a degree $(dim(S')-dim(S))$ morphism of graded
$\Q$-vector spaces
$$\pi_*:CH(X\times_{S'}Y,\Q)\rightarrow CH(X\times_{S}Y,\Q)$$
Let us show that these morphisms are compatible with the
composition of correspondences. Let $Z$ be another smooth
projective $S'$-variety, and let $\alpha$ and $\beta$ be elements
of $CH(X\times_{S'}Y,\Q)$, resp. $CH(Y\times_{S'}Z,\Q)$.

Consider the following diagram of Cartesian squares.
$$\begin{CD}
X\times_{S'}Y\times_{S'}Z@>i_2>> X\times_S (Y\times_{S'}Z) @>q_2
>>Y\times_{S'}Z
\\@Vi_1 VV @V\pi^X_2 VV @V\pi_2 VV
\\ (X\times_{S'}Y)\times_{S}Z@>\pi_1^Z>> X\times_S
Y\times_{S}Z@>p_2 >>Y\times_{S}Z
\\@Vq_1VV @Vp_1VV
\\ X\times_{S'}Y @>\pi_1 >> X\times_S Y
\end{CD}$$
and denote by $p_{XZ}$ the projection of $X\times_S Y\times_S Z$
on $X\times_S Z$, and by $p'_{XZ}$ the projection of $X\times_{S'}
Y\times_{S'} Z$ on $X\times_{S'} Z$. Finally, we denote by
$$\pi_*:CH(X\times_{S'}Z,\Q)\rightarrow CH(X\times_{S}Z,\Q)$$
the proper push-forward induced by $\pi:X\times_{S'}Z\rightarrow
X\times_S Z$.

 Using \cite[1.7]{FInt} (base change of the flat pull-back) and
\cite[8.3(c)]{FInt} (projection formula), we get
\begin{eqnarray*}
&&
\pi_*(p'_{XZ})_*[(i_1)^*(q_1)^*(\alpha)\,.\,(i_2)^*(q_2)^*(\beta)]
\\&=&(p_{XZ})_*(\pi_1^Z)_*(i_1)_*[(i_1)^*(q_1)^*(\alpha)\,.\,(i_2)^*(q_2)^*(\beta)]
\\&=& (p_{XZ})_*(\pi_1^Z)_*[(q_1)^*(\alpha)\,.\,(i_1)_*(i_2)^*(q_2)^*(\beta)]
\\&=&
(p_{XZ})_*(\pi_1^Z)_*[(q_1)^*(\alpha)\,.\,(\pi_1^Z)^*(\pi^X_2)_*(q_2)^*(\beta)]
\\&=& (p_{XZ})_*[(\pi_1^Z)_*(q_1)^*(\alpha)\,.\,(\pi^X_2)_*(q_2)^*(\beta)]
\\&=& (p_{XZ})_*[(p_1)^*(\pi_1)_*(\alpha)\,.\,(p_2)^*(\pi_2)_*(\beta)]
\end{eqnarray*}
in $CH(X\times_S\times Y\times Z,\Q)$. If $\alpha$ and $\beta$
belong to $Hom_{\mathcal{M}^o_+(S')}(M_{S'}(X),M_{S'}(Y))$, resp.
$Hom_{\mathcal{M}^o_+(S')}(M_{S'}(Y),M_{S'}(Z))$, then the first
of these expressions is the image of $\beta\circ \alpha$ under
$\pi_*$, while the last one is the composition of the images of
$\alpha$ and $\beta$ in
$Hom_{\mathcal{M}^o_+(S)}(M_{S}(X),M_{S}(Y))$, resp.
$Hom_{\mathcal{M}^o_+(S)}(M_{S}(Y),M_{S}(Z))$.

Hence, we obtain a $\Q$-linear forgetful functor
$$C\mathcal{V}^o_{S'}\rightarrow C\mathcal{V}^o_{S} $$
which passes to a $\Q$-linear forgetful functor
$$f_{*}:\mathcal{M}_{+}^o(S')\rightarrow \mathcal{M}_{+}^o(S)$$
by the universal property of the pseudo-abelian envelope.
\end{proof}
Observe that $f_*\circ f^*=M(S')\otimes (\,.\,)$.

We denote by $D^b(S,\Q_{\ell})$ the derived category of bounded
complexes of $\Q_{\ell}$-sheaves on $S$, in the sense of
\cite{adic}. The functor $\mathcal{V}_S\rightarrow
D^b(S,\Q_{\ell})$ mapping an object $h:X\rightarrow S$ to
$Rh_{*}(\Q_\ell)$ extends to a $\Q$-linear realization functor
$$real:\mathcal{M}_{+}^o(S)\rightarrow  D^b(S,\Q_{\ell})$$
By the relative K\"unneth isomorphism, it respects the tensor
structures: there is a canonical quasi-isomorphism
$$real\,(M_1\otimes M_2)\cong real\,(M_1)\otimes^L real\,(M_2)$$

\begin{definition}[Lefschetz motive]
Let $e:S\rightarrow \Pro^1_S$ be any section of the structural
morphism $\Pro^1_S\rightarrow S$. The cycle $\Pro^1_S\times_S
e(S)$ defines a projector $p$ on $M(\Pro^1_S)$ in the category
$\mathcal{M}^o_+(S)$, independent of the choice of $e$. The image
 of this projector is called the Lefschetz motive over $S$, and is
 denoted by $\LL_S$.
 \end{definition}

If $S'\rightarrow S$ is a morphism of smooth quasi-projective
varieties over $k$, then the base-change functor $\mot\rightarrow
\mathcal{M}^o_+(S')$ maps $\LL_S$ to $\LL_{S'}$. We will simply
write $\LL_S^d$ for
 $\LL_S^{\otimes d}$, for any integer $d\geq 0$. If
 $e_d:S\rightarrow\Pro^d_S$ is a section, we can define a
 projector $p_d=\Pro^d_S\times_S
e_d(S)$ on $M(\Pro^d_S)$ in the category $\mathcal{M}^o_+(S)$,
independent of the choice of $e_d$. Its image is isomorphic to
$\LL^d_S$, by the isomorphism $$\Pro^d _S\times_S e(S)^{\otimes
d}\in Hom_{\mathcal{M}^o_+(S)}((M(\Pro^d_S),p_d),\LL^d_S)$$

\begin{lemma}\label{twist}
For any pair of objects $X,Y$ in $\mathcal{V}_S$, with $X$ of pure
dimension $d$ over $S$, and for any pair of integers $r,s\geq 0$,
there is a canonical isomorphism of $\Q$-vector spaces
$$CH^{d+r-s}(X\times_S Y,\Q)\cong Hom_{\mathcal{M}^o_+(S)}(M(X)\otimes \LL^r_S,M(Y)\otimes \LL_S^{s})$$
where we put $CH^{j}(\,.\,)=0$ if $j<0$.
\end{lemma}
\begin{proof}
The functor $$(\,.\,)\otimes \LL_S:\mot\rightarrow \mot$$ is fully
faithful, by \cite[1.6]{deninger-murre}, so we may assume that
$r=0$ or $s=0$. For any integer $i>0$, we'll denote by
$f:X\times_S Y\times_S \Pro^i_S\rightarrow X\times_S Y$ the
projection morphism.

First, suppose $s=0$. An element of
$Hom_{\mathcal{M}^o_+(S)}(M(X)\otimes \LL^r_S,M(Y))$ is an element
of the form $\psi\circ (Id\otimes p_r)$, with $\psi$ in
$CH^{d+r}((\Pro^r_S\times_S X)\times_S Y,\Q)$. The computation in
\cite[p. 460]{manin} shows that the morphism
$$f^*:CH^{d+r}(X\times_S Y,\Q)\rightarrow CH^{d+r}((\Pro^r_S\times_S
X)\times_S Y,\Q)\circ (p_r\otimes Id)$$ is an isomorphism.

Now assume $r=0$. An element of
$Hom_{\mathcal{M}^o_+(S)}(M(X),M(Y)\otimes \LL^s_S)$ is an element
of the form $ (Id\otimes p_s)\circ \psi$, with $\psi$ in $CH^{d}(
X\times_S (Y\times \Pro^s_S),\Q)$. By an analogous computation,
the morphism
$$f_*: (Id\otimes p_s)\circ CH^{d}(X\times_S (Y\times_S \Pro^s_S),\Q)\rightarrow CH^{d-s}(X\times_S Y,\Q) $$ is an isomorphism.
\end{proof}

\begin{definition}[Chow group of a motive and cup product]
For any object $M$ in $\mathcal{M}^o_+(S)$, and any integer $d\geq
0$, we put
$$CH^d(M,\Q):=Hom_{\mathcal{M}^o_+(S)}(\LL^d,M)$$
If $\xi$ is an element of $CH^d(M,\Q)$, we denote by
$\xi\cup(\,.\,)$ the composition
$$\begin{CD}
M\otimes \LL^d_S@>Id\otimes \xi
>>M\otimes M@>\Delta_M >>M
\end{CD}$$
where $\Delta_{M}$ is the diagonal morphism (i.e. if $M$ is
$(M(Y),p)$ with $Y$ a smooth projective $S$-variety and $p$ a
projector on $M(Y)$, then $\Delta_M=p\circ M(\Delta_Y)\circ
(p\otimes p)$).
\end{definition}

\begin{lemma}\label{chowgr}
For any object $X$ in $\mathcal{V}_S$, for any pair of integers
$r,s\geq 0$, and for any object $N$ in $\mathcal{M}^o_+(S)$, there
is a canonical isomorphism
$$Hom_{\mathcal{M}^o_+(S)}(M(X)\otimes \LL^r_S,N\otimes \LL^s_S)
\cong \oplus_{X_i\in\mathcal{C}(X)}CH^{d(X_i/S)+r-s}(M(X_i)\otimes
N,\Q)$$ where $\mathcal{C}(X)$ denotes the set of connected
components of $X$, and $d(X_i/S)$ is the relative dimension of
$X_i$ over $S$.
\end{lemma}
\begin{proof}
We may suppose that $X$ has pure dimension $d$ over $S$. The
result holds if $N=M(Y)$ for some object $Y$ in $\mathcal{V}_S$,
by Lemma \ref{twist}. In general, $N$ is of the form $(M(Y),p)$,
with $p$ a projector on $M(Y)$. We have
\begin{eqnarray*}
CH^{d+r-s}(M(X)\otimes N,\Q)&\cong&(Id\otimes p)\circ
Hom_{\mathcal{M}^o_+(S)}(\LL^{d+r-s},M(X)\otimes M(Y))
\\&\cong&(p\otimes Id)\circ
Hom_{\mathcal{M}^o_+(S)}(M(X)\otimes \LL^r_S,M(Y)\otimes \LL^s_S)
\\&\cong &Hom_{\mathcal{M}^o_+(S)}(M(X)\otimes \LL_S^r,N\otimes \LL_S^s)
\end{eqnarray*}
\end{proof}

\subsection{Manin's identity principle}
\begin{lemma}\label{yoneda}
We denote by $Vec_{\Q}$ the category of $\Q$-vector spaces. Let
$N$ be any object of $\mathcal{M}^o_+(S)$, and consider the
functor
$$\omega_N:\mathcal{V}_S^{op}\rightarrow Vec_{\Q}$$
defined by $\omega_N(Y)=Hom_{\mot}(M(Y),N)$. The functor
$$\omega:\mot\rightarrow PreSh(\mathcal{V}_S):N\mapsto \omega_N$$
is fully faithful; here $PreSh(\mathcal{V}_S)$ denotes the
category of presheaves of $\Q$-vector spaces on $\mathcal{V}_S$.
\end{lemma}
\begin{proof}
By Yoneda's Lemma, the functor $$h:\mot\rightarrow
PreSh(\mot):N\mapsto h_N:= Hom_{\mot}(\,.\,,N)$$ is fully
faithful. The composition $\omega$ of $h$ with the restriction
functor
$$PreSh(\mot)\rightarrow PreSh(\mathcal{V}_S)$$ is still fully
faithful, since $h$ is recovered from $\omega$ by putting
$$h_{(M(Y),p)}(N):=\mathrm{Im}\{\,\omega_N(p):\omega_N(Y) \rightarrow \omega_N(Y)\,\}$$
\end{proof}
\begin{prop}[Manin's identity prinicple]\label{idpri}
Let $f,g:M\rightarrow N$ be morphisms in $\mot$.
\begin{itemize}
\item   The morphism $f$ is an isomorphism iff the morphism
$$\omega(f)(Y):CH^d(M(Y)\otimes M,\Q)\rightarrow CH^{d}(M(Y)\otimes N,\Q) $$ is an isomorphism for any integer $d\geq 0$, and any object $Y$ in $\mathcal{V}_S$ of
pure dimension $d$ over $S$. Moreover, $f=g$ iff
$\omega(f)(Y)=\omega(g)(Y)$ for all these pairs $(d,Y)$. \item A
sequence
$$\begin{CD} M@>a>> N @>b>> P
\end{CD}$$
in $\mot$ is split exact, iff $a$ admits a right inverse $a'$ and
the sequence
$$\begin{CD} CH^d(M(Y)\otimes M,\Q) @>\omega(a)(Y)>> CH^d(M(Y)\otimes N,\Q)
@>\omega(b)(Y)>> CH^d(M(Y)\otimes P,\Q)
\end{CD}$$ is a short exact sequence,
 for any integer $d\geq 0$, and any object $Y$ in $\mathcal{V}_S$ of
pure dimension $d$ over $S$.
\end{itemize}
\end{prop}
\begin{proof}
This follows from Lemma \ref{chowgr}, Lemma \ref{yoneda}, and the
proposition in \cite[p. 453]{manin}.
\end{proof}
\subsection{Relative Chow motive of a projective bundle}
Let $X$ be a smooth quasi-projective variety over $k$. In this
section, we compute the relative Chow motive of a projective
bundle $\pi:\Pro (\mathcal{E})\rightarrow X$, where $\mathcal{E}$
is a locally free sheaf of modules on $X$ of rank $r+1$. We show
that Manin's formula \cite[p. 457]{manin} in $\mathcal{M}^o_+(k)$
holds already in $\mathcal{M}^o_+(X)$. The proof closely follows
the arguments in \cite[2.4]{Scholl}.

\begin{prop}\label{bundle}
Let $X$ be a smooth quasi-projective variety over $k$, and let
$\mathcal{E}$ be a locally free sheaf of modules on $X$ of rank
$r+1$. Then
$$M(\Pro(\mathcal{E}))\cong \oplus_{i=0}^{r}\LL_X^{i}$$ in $\mathcal{M}^o_+(X)$.
\end{prop}
\begin{proof}
Consider the tautological bundle
$\mathcal{O}_{\Pro(\mathcal{E})}(1)$ on $\Pro(\mathcal{E})$, and
its divisor class $\xi:=c_1(\mathcal{O}_{\Pro(\mathcal{E})}(1))$
in $CH^1(\Pro(\mathcal{E}),\Q)$. By Lemma \ref{twist}, the element
$\xi^i\in CH^i(\Pro(\mathcal{E}))$ corresponds to a morphism
$\xi^i_*:\LL^i_X\rightarrow M(\Pro(\mathcal{E}))$ in
$\mathcal{M}^o_+(X)$, for any integer $i\geq 0$.
We
will show that the morphism
$$\psi:\oplus_i\xi^i_*: \oplus_i\LL^i_X\rightarrow
M(\Pro(\mathcal{E}))$$
is an isomorphism of relative motives over $X$.

By Manin's identity principle in Proposition \ref{idpri}, it
suffices to show that for any integer $d\geq 0$ and any object
$f:Y\rightarrow X$ in $\mathcal{V}_X$ of pure relative dimension
$d$,
$$\omega(\psi)(Y):\oplus_{i}CH^d(M(Y)\otimes \LL^i_X,\Q)\rightarrow CH^d(Y\times_X \Pro(\mathcal{E}),\Q) $$
is an isomorphism. However, $Y\times_X \Pro(\mathcal{E})\cong
\Pro(f^*\mathcal{E})$, and $f^*(\xi)$ is the divisor class $\xi'$
of $\mathcal{O}_{\Pro(f^*(\mathcal{E}))}$. Moreover,
$CH^d(M(Y)\otimes \LL^i_X,\Q)\cong CH^{d-i}(Y,\Q)$ by Lemma
\ref{chowgr}, and the morphism $CH^{d-i}(Y,\Q)\rightarrow
CH^d(\Pro(f^*\mathcal{E}),\Q)$ is intersection with $(\xi')^i$. So
it follows from Grothendieck's computation of the Chow group of a
projective bundle on a smooth quasi-projective variety
\cite[I.11]{semchev}, that $\omega(\psi)(Y)$ is an isomorphism.
\end{proof}
\begin{cor}\label{bundlecor}
Let $X$ be a smooth projective variety over $S$, and let
$\mathcal{E}$ be a locally free sheaf of modules on $X$ of rank
$r+1$. Denote by $\xi$ the divisor class
$c_1(\mathcal{O}_{\Pro(\mathcal{E})}(1))$ of the tautological line
bundle on $\Pro(\mathcal{E})$, and by
$h:\Pro(\mathcal{E})\rightarrow X$ the projection. Then
$$\sum_{i=0}^{r}\xi^i\cup(h\otimes Id)^*: \oplus_{i=0}^{r}(M(X)\otimes\LL_S^{i})\rightarrow M(\Pro(\mathcal{E}))$$ is an isomorphism in $\mathcal{M}^o_+(S)$.
\end{cor}
\begin{proof}
This is the image of the isomorphism in Proposition \ref{bundle}
under the forgetful functor $\mathcal{M}^o_+(X)\rightarrow \mot$.
\end{proof}
\subsection{The blow-up complex}\label{blowup}
Let $S$ be, as before, a smooth irreducible quasi-projective
variety over $k$. Let $X$ be a smooth projective irreducible
variety over $S$, and let $Z$ be a closed irreducible subvariety
of $X$ of codimension $r$, smooth over $S$. Consider the blow-up
$h:X'\rightarrow X$ of $X$ at $Z$, and denote by $E$ its
exceptional variety. We denote by $i$ and $j$ the inclusions of
$Z$ in $X$, resp. $E$ in $X'$, and we denote by $h_E:E\rightarrow
Z$ the restriction of $h$.
$$\begin{CD}
E@>j>> X'
\\ @Vh_E VV @VV h V
\\ Z@>i>> X
\end{CD}$$
We denote by $N$ the normal bundle of $Z$ in $X$.

Let $\zeta\in CH^{r-1}(E,\Q)$ be the $(r-1)$-th Chern class of the
quotient bundle $h^*(N)/\mathcal{O}_E(-1)$ on $E$. We denote by
$$i_*\in Hom_{\mot}(M(Z)\otimes \LL^r_S,
M(X))=CH^{d(Z/S)+r}(Z\times_S X,\Q)$$ the morphism defined by the
transpose of the graph of $i$. The morphism $$j_*:M(E)\otimes
\LL_S\rightarrow M(X')$$ is defined analogously. We'll write
$(h_E)_*$ for the morphism in
$$Hom_{\mot}(M(E)\otimes \LL_S, M(Z)\otimes
\LL^r_S)=CH^{d(E/S)-(r-1)}(E\times_S Z,\Q)$$ defined by the
transpose of the graph of $h_E$. With some abuse of notation,
we'll also write $\zeta\cup$ for the morphism
$M(E)\otimes\LL^r_S\rightarrow M(E)\otimes \LL_S$ obtained by
twisting $\zeta\cup$ (i.e. tensoring with the identity on
$\LL_S$).

\begin{lemma}\label{monoidal}
  Let $\xi$ be the
divisor class of the tautological bundle $\mathcal{O}_E(1)$ on
$E$. The morphism
$$\varphi:M(X)\oplus(\oplus_{i=1}^{r-1}M(Z)\otimes \LL_S^i) \rightarrow M(X')$$
defined by
$$\varphi(x,y_1,\ldots,y_{r-1})=h^*(x)+\sum_{i=1}^{r-1}j_*(\xi^{i-1}\cup (h_E\otimes Id)^*(y_i))$$
is an isomorphism in $\mot$.
\end{lemma}
\begin{proof}
By Manin's identity principle in Proposition \ref{idpri}, it
suffices to show that, for any integer $d\geq 0$, and any smooth
projective $S$-variety of pure relative dimension $d$,
$$\omega(\varphi)(Y):CH^d(Y\times_S X,\Q)\oplus(\oplus_{i=1}^{r-1}CH^{d-i}(Y\times_S Z,\Q)) \rightarrow
CH^d(Y\times_S X',\Q)$$ is an isomorphism.

Since blow-up commutes with flat base change, $h^Y:Y\times_S
X'\rightarrow Y\times_S X$ is the blow-up with center $Y\times_S
Z$ and exceptional divisor $Y\times_S E$. Moreover, the divisor
class of $\mathcal{O}_{Y\times_S E}(1)$ is the pull-back of the
divisor class of $\mathcal{O}_E(1)$. Hence, it suffices to show
that the morphism
\begin{eqnarray*}CH^d(
X,\Q)\oplus(\oplus_{i=1}^{r-1}CH^{d-i}(Z,\Q)) \rightarrow CH^d(
X',\Q)\,:\\(x,y_1,\ldots,y_{r-1})\mapsto
h^*(x)+\sum_{i=1}^{r-1}j_*(\xi^{i-1}\cup
h_E^*(y_i))\end{eqnarray*}
 is an isomorphism for each $d$. This is done in \cite[0.1.3]{beauv}.
\end{proof}

\begin{theorem}[Blow-up complex]\label{split}
Let $X$ be a smooth projective irreducible variety over $S$, and
let $Z$ be a closed irreducible subvariety of $X$ of codimension
$r$, smooth over $S$. Consider the blow-up $h:X'\rightarrow X$ of
$X$ at $Z$, and denote by $E$ its exceptional variety. We denote
by $i$ and $j$ the inclusions of $Z$ in $X$, resp. $E$ in $X'$,
and we denote by $h_E:E\rightarrow Z$ the restriction of $h$. The
complex
$$\begin{CD}
0@>>> M(X)@>h^*+ i^* >> M(X')\oplus M(Z)@>j^*-h_E^*>> M(E)@>>>
0\end{CD}$$ splits. In particular,
$$[M(X)]+[M(E)]=[M(X')]+[M(Z)]$$ in $K_0(\mathcal{M}_+^o(S))$.

If $G$ is a finite group, $X$ carries a good $G$-action, and $Z$
is $G$-closed, then the splitting is $G$-equivariant.
\end{theorem}
\begin{proof}
The absolute case $S=\mathrm{Spec}\,k$ was proven in
\cite[5.1]{GuiNav}.

  Let $\xi$ be the
divisor class of the tautological bundle $\mathcal{O}_E(1)$ on
$E$. We've seen in Corollary \ref{bundlecor} that
$$\psi(y_0,\ldots,y_r)=\sum_{i=0}^{r-1}\xi^{i}\cup (h_E\otimes Id)^*(y_i)$$ defines an isomorphism
$$\psi:\oplus_{i=0}^{r-1}M(Z)\otimes \LL_S^i \rightarrow M(E)$$
in $\mot$.
Combining this expression
with the isomorphism $\varphi$ in Lemma \ref{monoidal}, and
putting $A:=\oplus_{i=1}^{r-1}M(Z)\otimes \LL_S^i$, we get an
isomorphism of complexes
$$\begin{CD}
M(X)@>\left(\begin{array}{c}1\\ 0\\i^*\end{array}\right)>> M(X)\oplus A\oplus M(Z)@>\left(%
\begin{array}{c}
  0 \ \quad  1 \ \quad  0 \\
  i^* \quad 0 \  -1 \\
\end{array}%
\right)>> A\oplus M(Z)
\\@VIdVV @V(\varphi,Id) VV @VV\psi V
\\M(X)@>h^*+ i^* >> M(X')\oplus M(Z)@>j^*-h_E^*>>
M(E)
\end{CD}$$
Commutativity of the squares follows from the equality
$j^*j_*=\xi\cup$ (see the proof of \cite[6.7(c)]{FInt}). It is
clear that the upper complex splits.

If $X$ carries a good $G$-action, and $Z$ is $G$-closed, then all
the morphisms used in the proof are $G$-equivariant, and the
splitting is $G$-equivariant.
\end{proof}

\begin{remark}
In the classical setting, where $k=\C$, $S=\mathrm{Spec}\,\C$, and
where we use singular homology instead of motives, the exact
blow-up sequence arises as follows \cite[p.605]{GH}: take tubular
neighbourhoods $U$ and $V$ of $Z$ and $E$ in $X$, resp. $X'$. We
have $H_{*}(U)=H_{*}(Z)$, $H_{*}(V)=H_{*}(E)$,
$H_{*}(U-Y)=H_{*}(V-E)$, and $H_*(X'-E)=H_*(X-Z)$. Hence,
combining the Mayer-Vietoris sequences for $X=U\cup (X-Z)$ and
$X'=V\cup (X' - E)$, with the Barrat-Whitehead Lemma
\cite[17.4]{Algtop}, and the fact that $H_{*}(X')\rightarrow
H_{*}(X)$ is surjective, we obtain a short exact sequence
$$0\rightarrow H_{*}(E)\rightarrow H_{*}(Z)\oplus H_{*}(X')\rightarrow H_{*}(X)\rightarrow 0$$
\end{remark}

\subsection{Direct limits of categories}
Let $I$ be a directed set, and let $(\mathcal{C}_i)_{i\in I}$ be a
direct system of categories, with transition functors $F_{j,i}$
for $j\geq i$ in $I$.

\begin{definition}
We define the direct limit $\mathcal{C}$ of $(\mathcal{C}_i)_{i\in
I}$ as follows:
\begin{itemize}
\item Put $O:=\sqcup_{i\in I}Ob(\mathcal{C}_i)$. On this class, we
define an equivalence relation $\sim$ as follows: for any two
objects $(i,\alpha)$ and $(j,\beta)$, with $\alpha$, $\beta$
objects of $\mathcal{C}_{i}$, resp. $\mathcal{C}_{j}$,
$(i,\alpha)\sim (j,\beta)$ iff there exists an index $k\geq i,j$
in $I$ such that $F_{k,i}(\alpha)=F_{k,j}(\beta)$. We put
$Ob(\mathcal{C}):=O/\sim$. \item Consider two objects $A$ and $B$
in $Ob(\mathcal{C})$, represented by $(i,\alpha)$, resp.
$(j,\beta)$ in $O$. Put
$$H((i,\alpha),(j,\beta)):=\{ (k,f)\,|\,k\in
I,\,\mbox{with}\,k\geq i,j,\,\mbox{and}\,f\in
Hom_{\mathcal{C}_k}(F_{k,i}(\alpha),F_{k,j}(\beta))\}$$ On this
class, we define an equivalence relation $\sim$ as follows:
$(k,f)\sim (k',f')$ iff there exists an element $\ell$ in $I$,
with $\ell\geq k$ and $\ell\geq k'$, such that
$F_{\ell,k}(f)=F_{\ell,k'}(f')$. Observe that this equivalence
relation is compatible with the composition of morphisms. We put
$$Hom_{\mathcal{C}}(A,B):=H((i,\alpha),(j,\beta))/\sim $$ This
definition does not depend on the choice of the objects
$(i,\alpha)$ and $(j,\beta)$ in $O$ representing $A$ and $B$.
\end{itemize}
For each $i\in I$, there is a natural functor
$F_i:\mathcal{C}_i\rightarrow \mathcal{C}$, mapping an object
$\alpha$ to the class of $(i,\alpha)$, and mapping a morphism $f$
to $(i,f)$.
\end{definition}

If all the $\mathcal{C}_i$ are small categories, then
$\mathcal{C}$ is a direct limit in the category of small
categories. In any case, it has the following universal property:
$F_i=F_j\circ F_{j,i}$ for any pair $j\geq i$ in $I$, and
 if $\mathcal{D}$ is any category, and
$G_i:\mathcal{C}_i\rightarrow \mathcal{D}$ is a system of functors
such that $G_i=G_j\circ F_{j,i}$ for any pair $j\geq i$ in $I$,
then there exists a unique functor $G:\mathcal{C}\rightarrow
\mathcal{D}$ satisfying $G_i=G\circ F_i$ for any $i\in I$.

If the categories $\mathcal{C}_i$ are $\Q$-linear with tensor
structure, and if the transition functors respect $\Q$-linearity
and tensor structure, then $\mathcal{C}$ is a $\Q$-linear category
with tensor structure in a natural way, and the natural morphisms
$F_i$ respect $\Q$-linearity and tensor structure.

If $(\mathcal{C}_i)_{i\in I}$ is a direct system of pseudo-abelian
categories $\mathcal{C}_i$, with additive transition functors,
then
 $\mathcal{C}$ is pseudo-abelian, since additive functors preserve direct sums.

\subsection{Constructible Chow motives}\label{conchow}
Now, let $S$ be any $k$-variety. An admissible stratification is a
finite stratification $\mathscr{S}=\{S_1,\ldots,S_p\}$ of $S$ into
smooth, irreducible, quasi-projective, locally closed subvarieties
$S_i$. The admissible stratifications of $S$ form a directed set
(where $\mathscr{S}\leq \mathscr{S}'$ iff $\mathscr{S}'$ is a
refinement of $\mathscr{S}$). We define
$\mathcal{M}_{+}^o(S,\mathscr{S})$ as the direct product of
categories $\prod_{i=1}^{p}\mathcal{M}_+^o(S_i)$. Base change
induces restriction functors
$$F_{\mathscr{S}',\mathscr{S}}:\mathcal{M}^o_+(S,\mathscr{S})\rightarrow
\mathcal{M}^o_+(S,\mathscr{S}')$$ where $\mathscr{S}'$ is an
admissible refinement of $\mathscr{S}$. We obtain an inductive
system of $\Q$-linear categories with tensor structure, indexed by
the admissible stratifications $\mathscr{S}$ of $S$.

\begin{definition}\label{directlim}
We define the category $CMot_S$ of constructible Chow motives over
$S$, as the direct limit of the direct system
$(\mathcal{M}_{+}^o(S,\mathscr{S}))_{\mathscr{S}}$. The category
$CMot_S$ is $\Q$-linear, with tensor structure, and
pseudo-abelian. There are natural $\Q$-linear tensor functors
$$F_{\mathscr{S}}:\mathcal{M}_{+}^o(S,\mathscr{S})\rightarrow CMot_S$$

By Definition \ref{groth}, we can associate a Grothendieck ring
$K_0(CMot_S)$ to the additive tensor category $CMot_S$.
\end{definition}

\begin{lemma}\label{grothlim}
$$K_0(CMot_S)\cong\lim_{\stackrel{\longrightarrow}{\mathscr{S}}}\prod_{S_i\in \mathscr{S}}K_0(\mathcal{M}^o_+(S_i))$$
where $\mathscr{S}$ runs over the admissible stratifications of
$S$.
\end{lemma}
\begin{proof}
The canonical functors $F_{\mathscr{S}}$ induce ring morphisms
$$K_0(F_{\mathscr{S}}):\prod_{S_i\in \mathscr{S}}K_0(\mathcal{M}^o_+(S_i))\cong
K_0(\prod_{S_i\in \mathscr{S}}\mathcal{M}^o_+(S_i))\rightarrow
K_0(CMot_S)$$ and hence a ring morphism
$$\psi:\lim_{\stackrel{\longrightarrow}{\mathscr{S}}}\prod_{S_i\in \mathscr{S}}K_0(\mathcal{M}^o_+(S_i))\rightarrow K_0(CMot_S) $$
It is clear that $\psi$ is surjective, so let us prove
injectivity. Let $\mathscr{S}$ be an admissible stratification of
$S$. Let $\alpha$ be an element of
$K_0(\mathcal{M}^o_+(S,\mathscr{S}))$, and suppose that $\alpha$
maps to zero under $K_0(F_{\mathscr{S}})$. This means that there
exist objects $A$ and $B$ in $\mathcal{M}^o_+(S,\mathscr{S})$, and
an object $C$ in $CMot_S$, such that $\alpha+ [B]=[A]$, and
$F_{\mathscr{S}}(A)\oplus C$ and $F_{\mathscr{S}}(B)\oplus C$ are
isomorphic. By definition of the direct limit $CMot_S$, this
implies that there exists a refinement $\mathscr{S}'$ of
$\mathscr{S}$, and an object $\mathcal{C}'$ in
$\mathcal{M}^o_+(\mathscr{S}')$, such that
$F_{\mathscr{S}',\mathscr{S}}(A)\oplus C$ and
$F_{\mathscr{S}',\mathscr{S}}(B)\oplus C$ are isomorphic. Hence,
$\alpha$ maps to zero in
$$\lim_{\stackrel{\longrightarrow}{\mathscr{S}}}\prod_{S_i\in
\mathscr{S}}K_0(\mathcal{M}^o_+(S_i))$$
\end{proof}
\begin{cor}\label{realiz}
The realization functors $\mathcal{M}^o_+(S)\rightarrow
D^b(S,\Q_{\ell})$ for smooth quasi-projective $k$-varieties $S$,
induce a realization morphism of Grothendieck rings
$$K_0(CMot_S)\rightarrow  K_0(D^b(S,\Q_{\ell}))$$ for any
$k$-variety $S$. Here we take the Grothendieck ring of
$D^b(S,\Q_{\ell})$ as a triangulated category, i.e. if
$A\rightarrow B\rightarrow C\rightarrow A[1]$ is a distinguished
triangle in $D^b(S,\Q_{\ell})$, then $[B]=[A]+[C]$ in
$K_0(D^b(S,\Q_{\ell}))$.
\end{cor}
\begin{proof}
We define, for any admissible stratification $\mathscr{S}$ of $S$.
$$\psi_{\mathscr{S}}:K_0(\mathcal{M}^o_+(S,\mathscr{S}))\rightarrow
K_0(D^b(S,\Q_{\ell}))$$ as follows: if $S_i$ is a stratum of
$\mathscr{S}$, and $j:S_i\rightarrow S$ is the inclusion, then
$$K_0(\mathcal{M}^o_+(S_i))\rightarrow
K_0(D^b(S,\Q_{\ell}))$$ is the composition of
$$K_0(real):K_0(\mathcal{M}^o_+(S_i))\rightarrow
K_0(D^b(S_i,\Q_{\ell}))$$ with
$$K_0(j_{!}):K_0(D^b(S_i,\Q_{\ell}))\rightarrow
K_0(D^b(S,\Q_{\ell}))$$ These morphisms $\psi_{\mathscr{S}}$ form
a direct system and pass to a limit morphism
$$K_0(CMot_S)\rightarrow  K_0(D^b(S,\Q_{\ell}))$$
since for any closed subvariety $V$ of $S$, and any bounded
complex of $\Q_{\ell}$-sheaves $\mathcal{F}$ on $S$, we get a
distinguished triangle
$$
u_{!}u^*\mathcal{F}\rightarrow\mathcal{F}\rightarrow
v_{*}v^{*}\mathcal{F}\rightarrow u_{!}u^*\mathcal{F}[1]
$$
in $D^b(S,\Q_{\ell})$, where $U$ denotes the complement of $V$ in
$S$, and $u$ and $v$ are the inclusions of $U$, resp. $V$ in $S$.
\end{proof}

If $f:S'\rightarrow S$ is a morphism of smooth, quasi-projective
$k$-varieties, then the base-change functor $\mot\rightarrow
\mathcal{M}^o_+(S')$ induces a base-change functor
$CMot_S\rightarrow CMot_{S'}$. If $S'$ is a smooth projective
$S$-variety, the forgetful functor from Lemma \ref{forget} will in
general \textit{not} pass to a forgetful functor on constructible
Chow motives (due to the stratifications we allowed on the base).
\subsection{Extension of cohomological functors}\label{extension}
Let $G$ be a fixed finite group. Let $S$ be a variety over a field
$k$ of characteristic zero.

Theorem 2.2.2 in \cite{GuiNav} gives a very useful criterion to
extend cohomological functors defined on smooth and projective
varieties over $k$. We will extend this result to the relative
case, working over the base $S$.

Their cohomological theories take values in so-called
cohomological descent categories \cite[1.7.1]{GuiNav}. These are
tuples $(D,E,s,\lambda,\mu)$ which capture the essential
properties of the category of complexes over an abelian category,
and the class of quasi-isomorphisms. For our purposes, it suffices
to recall that one can associate a cohomological descent category
to any additive category $\Aa$, with underlying category
$D=\mathcal{C}^b(\Aa)$, where $s$ is the functor that associates
to a codiagram of complexes its total complex in
$\mathcal{C}^b(\Aa)$, and where $E$ is the class of homotopy
equivalences, by \cite[1.7.7]{GuiNav}. For any descent category
$(D,E,s,\lambda,\mu)$, we denote by $Ho(D)$ the localization of
$D$ w.r.t. the class of morphisms $E$.

We denote by $(G,Var_S)$ the category of varieties over $S$ with
good $G$-action. Let $(G,\mathcal{V}_{S})$ be the full subcategory
of $(G,Var_S)$ whose objects are the smooth and projective
varieties over $S$. We will also consider
the category $(G,Var_{S,c})$ of $S$-varieties with good
$G$-action, with proper morphisms. If $(G,\mathcal{C}_S)$ is any
of these categories, we define $(G,\mathcal{\overline{C}}_S)$ as
the inductive limit of
$(G,\mathcal{C}_{\mathscr{S}}):=\prod_{S_i\in
\mathscr{S}}(G,\mathcal{C}_{S_i})$ over the admissible
stratifications $\mathscr{S}=\{S_i\}_i$ of $S$ (see Section
\ref{conchow}). Beware: if $S_i$ is a locally closed subvariety of
$S$, and $X$ is an $S$-variety, the restriction of $X$ to $S_i$ is
given by $(X\times_S S_i)_{red}$ with its reduced structure.

\begin{definition}[Acyclic diagrams]
 An acyclic diagram in $(G,{Var}_{\mathscr{S}})$, for some admissible
stratification $\mathscr{S}=\{S_i\}$ of $S$, is a set of
$G$-equivariant Cartesian diagrams of the type

$$\begin{CD} \label{acyc}\tilde{Y}_i @>>> \tilde{X}_i
\\ @V(1)\hspace{140pt}VV @VV\pi V
\\Y_i @>>> X_i
\end{CD}$$

\noindent where $X_i,\, Y_i$ are $S_i$-varieties with good
$G$-action, the horizontal arrows are closed immersions, the
vertical arrows are proper $S_i$-morphisms, and $\pi$ induces an
isomorphism $\tilde{X}_i\setminus\tilde{Y}_i\cong X_i\setminus
Y_i$.

An acyclic diagram in $\overline{Var}_S$ is the image of an
acyclic diagram in $(G,Var_{\mathscr{S}})$ under the natural
functor $$F_{\mathscr{S}}:(G,Var_{\mathscr{S}})\rightarrow
(G,\overline{Var}_S)$$ for some admissible stratification
$\mathscr{S}$ of $S$.

An acyclic diagram in
 $(G,\overline{Var}_{S,c})$ is a diagram in
 $(G,\overline{Var}_{S,c})$ which is acyclic in $(G,\overline{Var}_{S})$.

 An elementary acyclic diagram in
$(G,\mathcal{\overline{V}}_S)$ is a blow-up diagram in
$(G,\mathcal{\overline{V}}_S)$, i.e. an acyclic diagram in
$(G,\overline{Var}_S)$ where all objects belong to
$(G,\mathcal{\overline{V}}_S)$, and $\pi$ is the blow-up of $X_i$
at the $G$-closed center $Y_i$.

An (elementary) acyclic morphism in any of the above categories,
is a morphism that can be realized as the right vertical arrow in
an (elementary) acyclic diagram.
\end{definition}

The following Theorem is a relative version of \cite{GuiNav},
Theorem 2.2.2 and \cite{DelAz}, Theorem 3.1. Let
$(D,E,s,\mu,\lambda)$ be a cohomological descent category.

\begin{theorem}[Extension Theorem]\label{descent}
Given a contravariant functor $F$ from
$(G,\mathcal{\overline{V}}_S)$ to a cohomological descent category
$(D,E,s,\mu,\lambda)$ such that
\begin{enumerate}
\renewcommand{\labelenumi}{({F}\theenumi)}
\item $F(\emptyset)=0$, \item the natural morphism $F(X\sqcup
Y)\rightarrow F(X)\times F(Y)$ is an isomorphism, \item for any
elementary acyclic diagram $X_{\bullet}$, the object
$sF(X_{\bullet})$ is acyclic in $Ho(D)$,
\end{enumerate}
there exists an essentially unique extension of $F$ to a functor
$$F_c:(G,\overline{Var}_{S,c})\rightarrow Ho(D)$$
such that the following descent properties hold:

(D) for any acyclic diagram $X_{\bullet}$ in
$(G,\overline{Var}_{S,c})$, the object $sF(X_{\bullet})$ is
acyclic.

(E) for any $S$-variety $X$ with good $G$-action, and any
$G$-closed, closed subvariety $Y$, we have a natural isomorphism
$$F_c(X\setminus Y)\cong s(F_c(Y\hookrightarrow X)) $$

 \noindent Moreover, if $H:\mathcal{\overline{V}}_S\rightarrow D$ is
another functor satisfying $(F_1),\,(F_2)$ and $(F_3)$, $H'$ is an
extension of $H$ satisfying $(D)$ and $(E)$, and
$\tau:F\rightarrow H$ is any natural transformation, then $\tau$
extends uniquely to a natural transformation $\tau:F'\rightarrow
H'$. This extension $\tau'$ is an isomorphism if $\tau$ is.
\end{theorem}

\begin{proof}
 As
explained in \cite{DelAz}, Proof of Theorem 2.2, and in the
appendix of \cite{DL5}, it suffices to use Proposition
\ref{resolution}, Proposition \ref{chow}, and Proposition
\ref{compact}, combined with
 the fact that all the objects in $\overline{Var}_{S,c}$ of
dimension zero are in $\overline{\mathcal{V}}_{S}$ (we de not
claim that any variety of relative dimension $0$ over $S$ is
smooth over $S$, but merely that this is true after admissible
stratification of the base $S$).

Now one can use the methods from \cite{GuiNav} (in particular, the
theory of cubical hyperresolutions) to prove Theorem
\ref{descent}.
\end{proof}

\subsection{Associating a constructible motive to a family of varieties with group
action}\label{groupmot} Let $G$ be a fixed finite group. Let $S$
be a variety over a field $k$ of characteristic zero.

For any additive category $\mathcal{A}$, we denote by
$(G,\mathcal{A})$ the additive category of functors $G\rightarrow
\Aa$, where we view $G$ as a category with one object, and
automorphism group $G$. We denote by
$Ho(G,\mathcal{C}^b(\mathcal{A}))=Ho\mathcal{C}^b(G,\mathcal{A})$
the homotopy category associated to $(G,\mathcal{A})$, i.e. the
category of bounded complexes over $(G,\mathcal{A})$ localized
w.r.t. homotopy equivalences.

If $S$ is a smooth, quasi-projective variety over $k$, we can
define a contravariant functor
$$M_G:(G,\mathcal{V}_S)\rightarrow (G,\mathcal{M}^o_+(S)) $$
mapping a smooth and projective $S$-variety $X$ with good
$G$-action to its motive $M(X)$, that we endow with a left
$G$-action as follows: an element $g$ of $G$ acts by its graph
$[g]$ in $Aut_{\mathcal{M}^o_+(S)}(M(X))$.

  For an arbitrary $k$-variety $S$, and any admissible stratification $\mathscr{S}$ of $S$, we have a
canonical functor
$$M_G:(G,\mathcal{V}_{\mathscr{S}})\rightarrow (G,CMot_S)$$ induced
by the functors $M:\mathcal{V}_{S_i}\rightarrow
\mathcal{M}^o_+(S_i)$ for $S_i\in \mathscr{S}$. These induce a
canonical functor
$$M_G:(G,\overline{\mathcal{V}}_S)\rightarrow (G,CMot_S)$$

\begin{theorem}
The functor
$$M_G:(G,\overline{\mathcal{V}}_S)\rightarrow (G,CMot_S)$$
has an essentially unique extension to a functor
$$M_{G,c}:(G,\overline{Var}_{S,c})\rightarrow Ho(G,\mathcal{C}^b(CMot_S))$$
satisfying properties $(D)$ and $(E)$ in Theorem \ref{descent}.
\end{theorem}
\begin{proof}
We only have to check properties $(F1),(F2),(F3)$ in Theorem
\ref{descent}. While $(F1)$ and $(F2)$ are obvious, $(F3)$ follows
from Theorem \ref{split}.
\end{proof}

\begin{cor}
There exists a functor
$$M_{G,c}:(G,Var_{S,c})\rightarrow Ho(G,\mathcal{C}^b(CMot_S))$$
that extends the functor $M_G:(G,\mathcal{V}_S)\rightarrow
(G,\mathcal{M}^o_+(S))$ and satisfies properties $(D)$ and $(E)$.

In particular, there exists a unique morphism of Grothendieck
rings
$$\chi_{G,c}:K^G_0(Var_S)\rightarrow K_0(G,CMot_S)$$ mapping the class of
an object $X$ of $\mathcal{V}_T$, with $T$ a smooth, irreducible,
quasi-projective locally closed subset of $S$, to the class of
$M_G(X)$ in $K_0(G,CMot_S)$.
\end{cor}
\begin{proof}
Use the canonical functor $F_{\{S\}}:(G,Var_{S,c})\rightarrow
(G,\overline{Var}_{S,c})$.

The existence of the morphism of Grothendieck rings $\chi_{c,G}$
follows from the application
$Ob\,Ho(\mathcal{C}^b(\Aa))\rightarrow K_0(\Aa)$ constructed in
\cite[5.4]{GuiNav} for any pseudo-abelian category $\Aa$.
\end{proof}

\begin{cor}\label{motives2}
There exists a functor
$$M_{c}:Var_{S,c}\rightarrow Ho(\mathcal{C}^b(CMot_S))$$
that extends the functor $M:\mathcal{V}_S\rightarrow
\mathcal{M}^o_+(S)$ and satisfies properties $(D)$ and $(E)$ (for
$G=\{e\}$).

In particular, there exists a unique morphism of Grothendieck
rings
$$\chi_{c}:K_0(Var_S)\rightarrow K_0(CMot_S)$$ mapping the class of
a smooth, projective $T$-variety $X$, with $T$ a smooth,
irreducible, quasi-projective locally closed subset of $S$, to the
class of $M(X)$ in $K_0(CMot_S)$.
\end{cor}

\begin{remark}\label{bittner}
In the Appendix to \cite{Bitt}, Bittner used the Weak
Factorization Theorem \cite{Weak} to prove that $K_0^{G}(Var_S)$
can be represented as follows: the set of generators consists of
the isomorphism classes $[X]$ of $S$-varieties $X$ with good
$G$-action, which are projective and smooth over their image in
$S$, and such that $G$ acts transitively on the connected
components of $X$. These generators are subject to the following
relations:
\begin{itemize}
\item $[\emptyset]=0$, \item$[Bl_{Y}X]-[E]=[X]-[Y]$, where $Y$ is
a closed $G$-invariant subvariety of $X$, smooth over its image in
$S$, which coincides with the image of X in $S$; and $Bl_{Y}X$ is
the blow-up of $X$ along $Y$, with exceptional divisor $E$, \item
$[X]=[X_T]+[X_{S\setminus T}]$, where $T$ is a closed subvariety
of $S$.
\end{itemize}
The existence and uniqueness of $\chi_{G,c}$ follow immediately
from this presentation.
\end{remark}

\begin{prop}
The functor $M_{G,c}$ commutes with base change, i.e. if
$S'\rightarrow S$ is a morphism of smooth quasi-projective
varieties, the diagram
$$\begin{CD}
(G,\overline{Var}_{S,c})@>M_{G,c}>>Ho(G,\mathcal{C}^b(CMot_S))
\\@VVV @VVV
\\(G,\overline{Var}_{S',c})@>M_{G,c}>>Ho(G,\mathcal{C}^b(CMot_{S'}))
\end{CD}$$
commutes (up to natural isomorphism). In particular, the diagram

$$\begin{CD}
K_0^G(Var_S)@>\chi_{G,c}>>K_0(G,CMot_S)
\\@VVV @VVV
\\K_0^G(Var_{S'})@>\chi_{G,c}>>K_0(G,CMot_{S'})
\end{CD}$$
 commutes.
\end{prop}

\begin{proof}
It is clear that both compositions coincide on
$(G,\overline{\mathcal{V}}_S)$. Moreover, the base change functor
$(G,Var_{S,c})\rightarrow (G,Var_{S',c})$ respects (elementary)
acyclic diagrams. So the diagram commutes by the uniqueness result
in Theorem \ref{descent}.
\end{proof}

\section{Character decomposition of constructible motives}\label{group} Throughout this section, our base scheme $S$
is a $k$-variety, with $k$ a field of characteristic zero. We will
generalize some results from \cite{DelAz} (where
$S=\mathrm{Spec}\,k$) to the relative, constructible setting, and
we will closely follow their arguments.

For any finite group $G$, and any $S$-variety $X$, we will denote
by $G\times X$ the constant finite group $S$-scheme associated to
$G$ and $X$. This means that $G\times X:=\sqcup_{g\in G}X$, and
$G$ acts on the left and on the right by permuting the indices. We
denote by $\mathbf{1}[G]$ the motive $M(G\times S)$ in $CMot_S$.
The left and right action of $G$ on $\mathbf{1}[G]$ are called the
left, resp. right regular representation of $G$. The object
$\mathbf{1}[G]$ is the motivic counterpart of the object $\Q[G]$
in the representation theory of $G$ over $\Q$.

If $M$ is an object of $CMot_S$ and $N$ is an object of
$(G,CMot_S)$, then we view $Hom_{CMot_S}(M,N)$ as a left
$\Q[G]$-module, by $(g,f)\mapsto g\circ f$.
\subsection{Restriction and induction}
Let $\psi:G\rightarrow G'$ be a morphism of finite groups. There
are obvious restriction functors, both denoted by $Res_{\psi}$,
$$
(G',\overline{Var}_{S,c})\rightarrow (G,\overline{Var}_{S,c}) $$
\vspace{-15pt}$$Ho(G',\mathcal{C}^b(CMot_S))\rightarrow
Ho(G,\mathcal{C}^b(CMot_S))$$

\begin{prop}\label{commuteres}
 The square
$$\begin{CD}
(G',\overline{Var}_{S,c}) @>M_{G',c}>>
Ho(G',\mathcal{C}^b(CMot_S))
\\ @VRes_{\psi}VV @VVRes_{\psi}V \\ (G,\overline{Var}_{S,c}) @>M_{G,c}>> Ho(G,\mathcal{C}^b(CMot_S))
\end{CD}$$
commutes (up to natural isomorphism).
\end{prop}
\begin{proof}
For objects in $(G',\overline{\mathcal{V}}_S)$, this is clear. By
Theorem \ref{descent}, the extension of
$$M_G\circ Res_{\psi}\cong Res_{\psi}\circ
M_{G'}:(G',\overline{\mathcal{V}}_S)\rightarrow
(G,\mathcal{C}^b(CMot_S))$$ to $(G',\overline{Var}_{S,c})$
satisfying $(D)$ and $(E)$ is essentially unique. Hence,
$M_{G,c}\circ Res_{\psi}\cong Res_{\psi}\circ M_{G',c}$.
\end{proof}

\begin{lemma}\label{reduc}
If $M$ is an object of $(G',CMot_S)$, and $N$ is an object of
$CMot_S$, then $Hom_{CMot_S}(N,M)$ is a left $\Q[G']$-module, and
$Hom_{CMot_S}(N,Res_{\psi}M)$ is the usual restriction
$Res_{\psi}Hom_{CMot_S}(N,M)$.
\end{lemma}
\begin{proof}
This is obvious.
\end{proof}

We construct a left adjoint $Ind_{\psi}$ for $Res_{\psi}$ as in
\cite{DelAz}, Section 4.

 If $X$ is an $S$-variety with good $G$-action, then
$G$ acts on $G'\times X$ by $g(g',x):=(g'\psi(g^{-1}),gx)$, and
$G'$ acts on $G'\times X$ by $g(g',x):=(gg',x)$. These actions
commute, and we put $Ind_{\psi}X:=(G'\times X)/G$ with its natural
$G'$-action. This construction commutes with stratification of the
base, and we obtain a functor
$$Ind_{\psi}:(G,\overline{Var}_{S,c})\rightarrow
(G',\overline{Var}_{S,c})$$

For any $\Q$-linear pseudo-abelian category $\Aa$, and any object
$M$ of $(G,\Aa)$, we can consider the projector
$\frac{1}{|G|}\sum_{g\in G}[g]$ on $M$ (where $[g]$ denotes the
image of $g$ in the endomorphism class of $M$). Its image is
denoted by $M^G$, and is called the $G$-invariant part of $M$. If
$N$ is an object of $\Aa$, then $G$ acts on the $\Q$-vector space
$Hom_{\Aa}(N,M)$, and $Hom_{\Aa}(N,M^G)$ coincides with the
subspace of $G$-invariants $(Hom_{\Aa}(N,M))^G$. We obtain an
additive functor $(.)^G:(G,\Aa)\rightarrow \Aa$.

If $M$ is an object of $(G,CMot_S)$, we let $G$ act on
$\mathbf{1}[G']\otimes M$ via the inverse of the right regular
representation of $G'$ and the $G$-action on $M$, and we let $G'$
act on $\mathbf{1}[G']\otimes M$ via the left regular
representation of $G'$ and the trivial action on $M$. These action
commute, and we obtain an object $Ind_{\psi}M:=(\1[G']\otimes
M)^G$ of $(G',CMot_S)$. We obtain additive functors, all denoted
by $Ind_{\psi}$,
$$
(G,CMot_S)\rightarrow (G',CMot_S) $$
\vspace{-15pt}$$(G,\mathcal{C}^b(CMot_S))\rightarrow
(G',\mathcal{C}^b(CMot_S))$$
\vspace{-15pt}$$Ho(G,\mathcal{C}^b(CMot_S))\rightarrow
Ho(G',\mathcal{C}^b(CMot_S))$$

If $H$ is a normal subgroup of $G$, and $\psi:G\rightarrow G/H$ is
the projection, then $Ind_{\psi}M$ is $M^H$ with its residual
$G/H$-action, and for any object $X$ in $(G,\overline{Var}_S)$,
$Ind_{\psi}X$ is $X/H$. In particular, for
$\psi:G\rightarrow\{e\}$, we obtain $M^G$ and $X/G$.

\begin{lemma}\label{induced}
If $M$ is an object of $(G,CMot_S)$, and $N$ is an object of
$CMot_S$, then $Hom_{CMot_S}(N,M)$ is a left $\Q[G]$-module, and
$Hom_{CMot_S}(N, Ind_{\psi}M)$ is the usual induced
$\Q[G']$-module $Ind_{\psi}Hom_{CMot_S}(N,M)$.
\end{lemma}
\begin{proof}
We see that
$$Hom_{CMot_S}(N,\1[G']\otimes M)\cong
\Q[G']\otimes_{\Q}Hom_{CMot_S}(N,M)$$ as a left $\Q[G']$-module.
Hence, it suffices to observe that, for any left $\Q[G]$-module
$V$,
$$(\Q[G']\otimes_{\Q}V)^G\cong \Q[G']\otimes_{\Q[G]}V$$ as a left
$\Q[G']$-module, where in the left hand side, $G$ acts on $\Q[G']$
via the inverse of the right regular representation.
\end{proof}

The following propositions follow immediately from the proofs of
their counterparts in \cite{DelAz}, Proposition 4.1-4.
\begin{prop}
For any object $X$ of $(G,\overline{Var}_S)$, and any object $Y$
of $(G',\overline{Var}_S)$, there is a natural bijection
$$Hom_{(G,\overline{Var}_S)}(X,Res_{\psi}Y)\cong Hom_{(G',\overline{Var}_S)}(Ind_{\psi}X,Y) $$
\end{prop}

\begin{prop}
The functors $Ind_{\psi}:(G,CMot_S)\rightarrow (G',CMot_S)$ and
$Ind_{\psi}:Ho(G,\mathcal{C}^b(CMot_S))\rightarrow
Ho(G',\mathcal{C}^b(CMot_S))$ are left adjoint to $Res_{\psi}$.
\end{prop}


\begin{prop}
For any object $M$ in $(G,\mathcal{C}^b(CMot_S))$, and any object
$N$ in $(G',\mathcal{C}^b(CMot_S))$, we have the projection
formula
$$Ind_{\psi}(Res_{\psi}N\otimes M)\cong N\otimes Ind_{\psi}M$$
\end{prop}

\begin{prop}\label{proj}
For any object $X$ of $(G,\overline{Var}_S)$, and any object $Y$
of $(G',\overline{Var}_S)$, we have the projection formula
$$Ind_{\psi}(Res_{\psi}Y\otimes X)\cong Y\otimes Ind_{\psi}X$$
\end{prop}

\subsection{The motive of a quotient variety}
\begin{theorem}\label{quot}
The diagram
$$\begin{CD} (G,\overline{Var}_{S,c}) @>M_{G,c}>>
Ho(G,\mathcal{C}^b(CMot_S))
\\ @VInd_{\psi}VV @VVInd_{\psi}V \\ (G',\overline{Var}_{S,c}) @>M_{G',c}>> Ho(G',\mathcal{C}^b(CMot_S))
\end{CD}$$
commutes.

 In other words, for any object $X$ in
$(G,\overline{Var}_{S})$, there exists a natural isomorphism
$M_{G',c}(Ind_{\psi}(X))\cong Ind_{\psi}(M_{G,c}(X))$.
\end{theorem}
\begin{proof}
In the absolute case $S=\mathrm{Spec}\,k$, this is the main result
of \cite{DelAz}. We will show that their arguments carry over the
the relative case. The result will appear as a consequence of the
uniqueness statement in Theorem \ref{descent}, if we can prove
that both paths in the diagram are isomorphic on
$(G,\overline{\mathcal{V}}_S)$, and satisfy $(D)$ and $(E)$ on
$(G,\overline{Var}_{S,c})$. This is done in Lemma \ref{isom} and
Lemma \ref{acycl} below.
\end{proof}
Let $X$ be any $S$-variety with good $G$-action. The quotient map
$$\pi:G'\times X\rightarrow (G'\times X)/G=Ind_\psi X$$
induces a morphism
$$\pi^*:M_{c,G}Ind_{\psi}X\rightarrow(\1[G']\otimes M_{G,c}(X))^G$$
and this defines a morphism of functors
$$\psi:M_{c,G}Ind_{\psi}\rightarrow Ind_{\psi}M_{c,G}$$
\begin{lemma}\label{isom}
The morphism $\psi$ is an isomorphism on
$(G,\overline{\mathcal{V}}_S)$.
\end{lemma}
\begin{proof}
 We will merely sketch the
arguments in \cite{DelAz}, to show that they carry over to our
setting. We may suppose that $G'=\{e\}$. Let $S$ be a smooth
quasi-projective variety over $k$. The first step is to construct
a category of effective Chow motives starting from $S$-varieties
of the form $X'=X/G$, with $X$ in $(G,\mathcal{V}_S)$. We denote
the full subcategory of $Var_S$ with these varieties as objects by
$\mathcal{V}_S'$. By \cite[17.4.10]{FInt}, the usual construction
of Chow motives still makes sense if we start from
$\mathcal{V}'_S$ instead of $\mathcal{V}_S$, and we obtain a
pseudo-abelian category $\mathcal{M}^o_+(S)'$. There is an obvious
fully faithful embedding
$$\Phi:\mot\rightarrow \mot '$$
and we show that it is an equivalence by establishing an
isomorphism between the motive of $X/G$ in $\mot'$, and the image
of $M(X)^G$ under $\Phi$. This is done as in \cite[1.2]{DelAz},
using Manin's identity principle and the fact that
$CH(X/G,\Q)\cong CH(X,\Q)^G$ (see \cite[1.7.6]{FInt}). We fix a
quasi-inverse functor for $\Phi$, and this yields a functor
$$M':\mathcal{V}_S'\rightarrow \mot :X/G\rightarrow M'(X/G)\cong M(X)^G$$
If $S$ is any $k$-variety, taking direct limits over admissible
stratifications of $S$ yields a functor
$$M':\overline{\mathcal{V}}_S'\rightarrow CMot_S$$

If we define elementary acyclic diagrams in
$\overline{\mathcal{V}}_S'$ as quotients of elementary acyclic
diagrams in $(G,\overline{\mathcal{V}}_S)$, then we can formulate
an extension principle for the category
$\overline{\mathcal{V}}_S'$ as in Theorem \ref{descent} (see
\cite[2.2]{DelAz} for the absolute case). In particular, $M'$ has
an essentially unique extension
$$M'_c:\overline{Var}_{S,c}\rightarrow Ho(\mathcal{C}^b(CMot_S))$$ satisfying $(D)$
and $(E)$. However, is is easily seen that $\psi$ defines an
isomorphism of functors $M\cong M'$ on $\overline{\mathcal{V}}_S$,
so by the uniqueness statement in Theorem \ref{descent}, $\psi$ is
also an isomorphism on $\overline{\mathcal{V}}_{S}'$. This
concludes the proof.
\end{proof}
\begin{lemma}\label{acycl}
The functors $$M_{c,G}Ind_{\psi}\mbox{\ and\ } Ind_{\psi}M_{c,G}:
(G,\overline{Var}_{S,c})\rightarrow Ho(G',\mathcal{C}^b(CMot_S))$$
satisfy properties $(D)$ and $(E)$ from Theorem \ref{descent}.
Hence,
$$\psi:M_{c,G}Ind_{\psi}\rightarrow Ind_{\psi}M_{c,G}$$ is an
isomorphism of functors on $\overline{Var}_{S,c}$.
\end{lemma}
\begin{proof}
We know that $M_{c,G}$ and $M_{c,G'}$ respect acyclic diagrams, by
construction. It is easy to see that $Ind_{\psi}$ respects acyclic
diagrams of varieties, and acyclic diagrams in the descent
categories (it commutes with $s$ since it is an additive functor).

The fact that $\psi$ is an isomorphism follows from Lemma
\ref{isom} and the uniqueness statement in Theorem \ref{descent}.
\end{proof}
\begin{cor}
For any $X$ in $(G,\overline{Var}_{S})$, $M_c(X/G)\cong
M_{G,c}(X)^G$.
\end{cor}

\subsection{Character decomposition}\label{chardec}
We denote by $C(G,\Q)$ the $\Q$-vector space of $\Q$-central
functions, i.e. the space of $\Q$-linear combinations of
characters of $\Q$-irreducible representations of $G$.
 Recall that a central function $\alpha:G\rightarrow \Q$
belongs to $C(G,\Q)$ iff $\alpha(x)=\alpha(x')$ for each pair
$\{x,x'\}$ of elements of $G$, for which the subgroups generated
by $x$, resp. $x'$, are conjugate in $G$. Artin proved that every
such $\alpha$ is a $\Q$-linear combination of characters of the
form $Ind^{G}_{H}1_H$, with $H$ a cyclic subgroup of $G$, and
$1_H$ the trivial character on $H$.

\begin{definition}
For any finite dimensional $\Q$-vector space $V$ with $G$-action,
we define an element $V$ of $(G,CMot_S)$ as follows: we write the
vector space $V$ as the image of $\Q[G]^n$ under a projector $p$,
and we define the object $V$ of $CMot_S$ as the image of the
corresponding projector on $\1[G]^{\oplus n}$.
\end{definition}

\begin{lemma}
For any pair of objects $M,N$ of $(G,CMot_S)$, and any finite
$\Q[G]$-module $V$, we have a natural isomorphism of
$\Q[G]$-modules
$$Hom_{CMot_S}(N,V\otimes M)\cong V\otimes_{\Q} Hom_{CMot_S}(N,M)$$
\end{lemma}
\begin{proof}
 If $V$ is the image of $\Q[G]^n$ under a projector $p$, we
have \begin{eqnarray*}
Hom_{CMot_S}(N,V\otimes M)&=&(p\otimes Id)\circ Hom_{CMot_S}(N,\1[G]^{\oplus n}\otimes M)\\
&=&p(\Q[G]^n)\otimes_{\Q}Hom_{CMot_S}(N,M) \\&=& V\otimes_{\Q}
Hom_{CMot_S}(N,M)
\end{eqnarray*}
\end{proof}

\begin{definition}
Let $\alpha$ be an effective character of $G$ over $\Q$, and let
$\rho_{\alpha}:G\rightarrow GL(V_{\alpha})$ be the corresponding
representation. We denote by $V_{\alpha}^{\vee}$ its dual
representation, as well as the associated object in $(G,CMot_S)$.
For any object $M$ of $(G,CMot_S)$, we put
$$M_{\alpha}:=(V_{\alpha}^{\vee}\otimes M)^G$$
in $CMot_S$. This defines additive functors, all denoted by
$(.)_{\alpha}$
$$
(G,CMot_S)\rightarrow CMot_S$$
\vspace{-15pt}$$(G,\mathcal{C}^b(CMot_S))\rightarrow
\mathcal{C}^b(CMot_S)$$
\vspace{-15pt}$$Ho(G,\mathcal{C}^b(CMot_S))\rightarrow
Ho(\mathcal{C}^b(CMot_S))$$
\end{definition}

\begin{lemma}\label{alpha}
For any $\Q$-irreducible character $\alpha$, and any pair of
objects $M,N$ of $(G,CMot_S)$, we have a natural isomorphism of
$\Q[G]$-modules
$$Hom_{CMot_S}(N,M_{\alpha})\cong
Hom_{\Q[G]}(V_\alpha,Hom_{CMot_S}(N,M))$$
\end{lemma}
\begin{proof}
\begin{eqnarray*}
Hom_{CMot_S}(N,M_{\alpha})&\cong&(V_{\alpha}^{\vee}\otimes_{\Q}Hom_{CMot_S}(N,M))^G
\\&\cong&Hom_{\Q}(V_{\alpha},Hom_{CMot_S}(N,M))^G
\\&\cong& Hom_{\Q[G]}(V_\alpha,Hom_{CMot_S}(N,M))
\end{eqnarray*}
\end{proof}

Suppose that $\alpha$ is $\Q$-irreducible, of degree $n_{\alpha}$.
We define an idempotent $p_{\alpha}$ in the group algebra $\Q[G]$
by
$$p_{\alpha}:=\frac{n_{\alpha}}{|G|<\alpha,\alpha>}\sum_{g\in
G}\alpha(g^{-1})[g]$$

\begin{prop}\label{definition}
For each $\alpha\in C(G,\Q)$, there exists a unique morphism of
abelian groups
$$\chi_{c,\alpha}:K_0^{G}(Var_S)\rightarrow K_0(CMot_S)\otimes
\Q\,$$ such that
\\(i) if $X$ belongs to $(G,\mathcal{V}_T)$, for some
smooth, irreducible, and quasi-projective locally closed subset
$T$ of $S$, and if $\alpha$ is a $\Q$-irreducible character on
$G$, then
$$\frac{n_{\alpha}}{<\alpha,\alpha>}\chi_{c,\alpha}([X])$$ is equal
to the class of the image of the projector $p_{\alpha}$ on the
object $M(X)$ of $\mathcal{M}^o_+(T)$.
\\(ii) the morphism $\chi_{c,\alpha}$ is $\Q$-linear in $\alpha$.
\end{prop}

\begin{proof}
We can construct $\chi_{c,\alpha}$ as follows: if $\alpha$ is
effective, we map the class of an $S$-variety $X$ with good
$G$-action in $K_0^G(Var_S)$, to the class of
$(M_{G,c}(X))_{\alpha}$ in $K_0(CMot_S)$. Property $(ii)$ is
clear, so let us prove $(i)$. For any object $N$ in $CMot_S$, we
have
$$Hom_{CMot_S}(N,\mathrm{Im}\,p_{\alpha}|_{M(X)})\cong p_{\alpha}\Q[G]\otimes_{\Q[G]}Hom(N,M(X)) $$
In view of Lemma \ref{alpha} and Yoneda's Lemma, any choice of an
isomorphism of $\Q$-vector spaces
$$\Q[G]\cong (V_{\alpha}^{\vee})^{\oplus n_{\alpha}/<\alpha,\alpha>} $$
induces an isomorphism
$$\mathrm{Im}_{p_{\alpha}}\cong (M(X)_{\alpha})^{\oplus n_{\alpha}/<\alpha,\alpha>}  $$
in $CMot_S$.

 Alternatively, we can take point $(i)$ as a
definition, and use Bittner's presentation of $K_0^G(Var_S)$ (see
Remark \ref{bittner} in Section \ref{groupmot}) and Theorem
\ref{split}.
\end{proof}

\begin{lemma}\label{prod}
Let $G_1$ and $G_2$ be finite groups, and let $\alpha_i$ be an
element of $C(G_i,\Q)$, for $i=1,2$. Consider the ring morphism
$$\Delta: K_0^{G_1}(Var_S)\otimes_{\Z} K_0^{G_2}(Var_S)\rightarrow K_0^{G_1\times
G_2}(Var_S)$$ obtained as follows: the projection of $G_1\times
G_2$ on $G_i$ induces a morphism
$\pi_i:K_0^{G_i}(Var_S)\rightarrow K_0^{G_1\times G_2}(Var_S)$.
The composition of
$$\pi_1\times \pi_2:K_0^{G_1}(Var_S)\times K_0^{G_2}(Var_S)\rightarrow K_0^{G_1\times
G_2}(Var_S)\times K_0^{G_1\times G_2}(Var_S)$$ with ring
multiplication in $K_0^{G_1\times G_2}(Var_S)$, is bilinear, and
induces the morphism $\Delta$. We have
$$\chi_{c,\alpha_1}\otimes \chi_{c,\alpha_2}=\chi_{c,\alpha_1.\alpha_2}\circ\Delta$$
\end{lemma}
\begin{proof}
Let $T$ be any smooth, quasi-projective, irreducible locally
closed subset of $S$. Let $X_i$ be a proper and smooth variety
over $T$, with good $G_i$-action, for $i=1,2$. We may assume
$\alpha_i$ is $\Q$-irreducible, for $i=1,2$. Let
$\rho_i:G_i\rightarrow GL(V_i)$ be the representation with
character $\alpha_i$. The external tensor product
$\rho_1\boxtimes\rho_2:G_1\times G_2\rightarrow
V_1\otimes_{\Q}V_2$ is an irreducible representation with
character $\alpha_1.\alpha_2$.

Hence, by Proposition \ref{definition}$(i)$, it suffices to show
that the projectors $Im\,p_{\alpha_1.\alpha_2}$ and
$Im\,p_{\alpha_1} \otimes Im\,p_{\alpha_2}$ on $M(X_1\times_T
X_2)$ have isomorphic images in $\mathcal{M}^o_+(T)$. This,
however, is clear.
\end{proof}


Now we prove two Frobenius reciprocity properties.

\begin{lemma}\label{welldef3}
Let $\psi:G\rightarrow G'$ be a morphism of finite groups. For any
$S$-variety $X$ with good $G$-action, and any $\alpha\in
C(G',\Q)$,
$$\chi_{c,\alpha}([Ind_{\psi}X])=\chi_{c,Res_{\psi}\alpha}([X])$$
in $K_0(CMot_S)\otimes\Q$.
\end{lemma}
\begin{proof}
We may assume that $\alpha$ is effective. By definition,
$$M_{G',c}(Ind_{\psi}X)_{\alpha}=(V_{\alpha}^{\vee}\otimes M_{G',c}(Ind_{\psi}X))^{G'}$$
By Proposition \ref{quot}, this object is isomorphic to
$$(V_{\alpha}^{\vee}\otimes Ind_{\psi}M_{G,c}(X))^{G'} $$
By the projection formula in Proposition \ref{proj}, we get
\begin{eqnarray*}
(V_{\alpha}^{\vee}\otimes Ind_{\psi}M_{G,c}(X))^{G'}&\cong&
(Ind_{\psi}(Res_{\psi}V_{\alpha}^{\vee}\otimes M_{G,c}(X)))^{G'}
\\&\cong& (V_{Res_{\psi}\alpha}^{\vee}\otimes M_{G,c}(X))^G
\\&\cong& M_{G,c}(X)_{Res_{\psi}\alpha}
\end{eqnarray*}
\end{proof}

\begin{cor}[Relative motive of a quotient variety]\label{quot2}
Let $X$ be an $S$-variety with good $G$-action, and let $H$ be a
normal subgroup of $G$. For every $\alpha$ in $C(G/H,\Q)$,
$$\chi_{c,\alpha}([X/H])=\chi_{c,\alpha\circ\rho}([X])\,$$
where $\rho$ is the projection $\rho:G\rightarrow G/H$.
\end{cor}
\begin{proof}
Apply Lemma \ref{welldef3} to the projection $\psi:G\rightarrow
G/H$.
\end{proof}

\begin{lemma}\label{welldef2}
Let $\psi:G\rightarrow G'$ be a morphism of finite groups. For any
$S$-variety $X$ with good $G'$-action, and any $\alpha\in
C(G,\Q)$,
$$\chi_{c,\alpha}([Res_{\psi}X])=\chi_{c,Ind_{\psi}\alpha}([X])$$
in $K_0(CMot_S)\otimes\Q$.
\end{lemma}
\begin{proof}
We may assume that $\alpha$ is effective. By definition,
\begin{eqnarray*}
M_{G',c}(X)_{Ind_{\psi}\alpha}&\cong&
(V_{Ind_{\psi}\alpha}^{\vee}\otimes M_{G',c}(X))^{G'}
\\&\cong& (Ind_{\psi}V_{\alpha}^{\vee}\otimes M_{G',c}(X))^{G'}
\end{eqnarray*}
By the projection formula in Proposition \ref{proj}, we get
\begin{eqnarray*}
(Ind_{\psi}V_{\alpha}^{\vee}\otimes M_{G',c}(X))^{G'}&\cong&
(Ind_{\psi}(V_{\alpha}^{\vee}\otimes Res_{\psi}M_{G',c}(X)))^{G'}
\\&\cong& (V_{\alpha}^{\vee}\otimes Res_{\psi}M_{G',c}(X))^{G}
\\&\cong& M_{G,c}(Res_{\psi}X)_{\alpha}
\end{eqnarray*}
where the last isomorphism follows from Proposition
\ref{commuteres}.
\end{proof}
\begin{corollary}\label{cor}
For every $S$-variety $X$ with good $G$-action,
$$\chi_{c}([X])=\sum_{\alpha}\frac{n_{\alpha}}{<\alpha,\alpha>}\chi_{c,\alpha}(X)\,$$
where we take the sum over all $\Q$-irreducible representations
$\alpha$.
\end{corollary}
\begin{proof}
The regular character $\chi_{reg}$ on $G$ is induced by the
trivial character on the trivial subgroup. Now use Lemma
\ref{welldef2}, and the fact that for each $\Q$-irreducible
character $\alpha$, its degree $n_{\alpha}$ equals its
multiplicity in $\chi_{reg}$ times $<\alpha,\alpha>$.
\end{proof}

\begin{definition}
We denote by $K_0^{mot}(Var_S)$ the image of the morphism
$$\chi_c:K_0(Var_S)\rightarrow K_0(CMot_S)$$ constructed in
Corollary \ref{motives2}.
\end{definition}

\begin{lemma}\label{artin}
For any $\Q$-central function $\alpha$ in $C(G,\Q)$, the image of
$\chi_{c,\alpha}$ is contained in $K_0^{mot}(Var_S)\otimes \Q$.
\end{lemma}
\begin{proof}
In view of Artin's result mentioned above, it suffices to prove
this lemma when $\alpha$ is of the form $Ind^{G}_H1_H$, where $H$
is a cyclic subgroup of $G$, and $1_H$ is the trivial character on
$H$. By Lemma \ref{welldef2}, we may assume that $G$ is cyclic,
and $G=H$. In this case, Lemma \ref{artin} is an immediate
consequence of Corollary \ref{quot2}.
\end{proof}

\begin{lemma}[Base Change]\label{basechange0}
Let $f:T\rightarrow S$ be a morphism of $k$-varieties. For any
$\alpha\in C(G,\Q)$, the base change square
$$\begin{CD}
K_0^G(Var_S)@>\chi_{c,\alpha}>>K_0^{mot}(Var_S)\otimes\Q
\\ @Vf^{*}VV @VVf^{*}V
\\ K_0^G(Var_T)@>\chi_{c,\alpha}>>K_0^{mot}(Var_T)\otimes\Q
\end{CD}$$ commutes.
\end{lemma}
\begin{proof}
This is clear from the definition.
\end{proof}
\section{A motivic incarnation for the relative theory of pseudo-finite
fields}\label{relative} Now we're ready to construct the relative
counterpart of the morphism $$\chi_{(c)}:K_{0}(PFF_k)\rightarrow
K_0^{mot}(Var_k)\otimes \Q$$ from the introduction. Throughout
this section, our base scheme $S$ is a variety over a field $k$ of
characteristic zero.

\begin{definition}\label{Theta}
Let $X$ be an affine normal irreducible variety over $S$, let $Y$
be a Galois cover of $X$, and let $Con$ be a conjugation domain of
the cover $Y/X$. Put $G=G(Y/X)$. We define a $\Q$-central function
$\alpha_{Con}$ in $C(G,\Q)$ as follows: $\alpha(g)=1$ if the
subgroup generated by $g$ belongs to $Con$, and $\alpha(g)=0$
else.

Let $\Theta$ be the set of quantifier-free Galois formulas
$\theta$ over $S$. We define a map
$$\chi_{(c)}:\Theta\rightarrow K_0^{mot}(Var_S)\otimes \Q$$ as
follows: for any quantifier-free Galois formula $\theta$ over $S$,
corresponding to a Galois stratification $<X,C_i/A_i,Con(A_i)>$,
we define $\chi_{(c)}(\theta)$ by $$\chi_{(c)}(\theta):=\sum_i
\chi_{c,\alpha_{Con(A_i)}}([C_i])$$
\end{definition}

\begin{lemma}\label{inflate}
Let $\theta$ and $\theta'$ be quantifier-free Galois formulas,
such that the Galois stratification corresponding to $\theta'$ is
obtained from the one corresponding to $\theta$ by inflation. Then
$$\chi_{(c)}(\theta)=\chi_{(c)}(\theta')$$
\end{lemma}
\begin{proof}
Let $X,Y$ be affine normal irreducible varieties over $S$, let
$Y/X$ be a Galois cover, and let $Con(Y/X)$ be a conjugation
domain for this cover. Put $G'=G(Y/X)$. Let $Z/X$ be a cover
dominating $Y/X$, put $G=G(Z/X)$, and denote by $Con(Z/X)$ the
conjugation domain obtained by inflation. Denote by
$\psi:G\rightarrow G'$ the projection.

We have $Y=Ind_{\psi}Z$, and hence, by Lemma \ref{welldef3}, it
suffices to observe that
$\alpha_{Con(Z/X)}=Res_{\psi}\alpha_{Con(Y/X)}$.
\end{proof}

\begin{lemma}\label{refine}
Let $\theta$ and $\theta'$ be quantifier-free Galois formulas,
such that the Galois stratification corresponding to $\theta'$ is
obtained from the one corresponding to $\theta$ by refinement.
Then
$$\chi_{(c)}(\theta)=\chi_{(c)}(\theta')$$
\end{lemma}

\begin{proof}
Let $X,Y$ be affine normal irreducible varieties over $S$, let
$Y/X$ be a Galois cover, and let $Con(Y/X)$ be a conjugation
domain for this cover. Put $G'=G(Y/X)$. Let $U$ be a normal
irreducible closed subvariety of $X$, and let $V$ be any connected
component of $Y_U=Y\times_X U$. Put $G=G(V/U)$, and let
$\psi:G\rightarrow G'$ be the inclusion. Let $Con(V/U)$ be the
conjugation domain induced from $Con(Y/X)$ by refinement.

Since $Y_U=Ind_{\psi}V$ and
$\alpha_{Con(V/U)}=Res_{\psi}\alpha_{Con(Y/X)}$, we see from Lemma
\ref{welldef3} that
$$M_{c,G}(Y_U)_{\alpha_{Con(Y/X)}}\cong
M_{c,G'}(V)_{\alpha_{Con(V/U)}}$$ The result now follows from
additivity of $M_{c,G}$.
\end{proof}

\begin{lemma}\label{welldef}
Let $m$ be a positive integer, let $\mathcal{A}$ and $\mathcal{B}$
be Galois stratifications of $\A^m_S$, and suppose that there
exists a Galois stratification $\mathcal{G}$ of $\A^m_S\times_S
\A^m_S$, such that, for each point $x$ of $S$,
 and each pseudo-finite field
extension $M$ of $k(x)$, the set $Z(\mathcal{G},x,M)$ is the graph
of a bijection between $Z(\mathcal{A},x,M)$ and
$Z(\mathcal{B},x,M)$. Let $\theta_A$,
 $\theta_B$, and $\theta_G$ be the quantifier-free Galois formulas corresponding to $\mathcal{A}$,
$\mathcal{B}$, resp. $\mathcal{G}$. Then
$$\chi_{(c)}(\theta_A)=\chi_{(c)}(\theta_B)=\chi_{(c)}(\theta_G)$$
\end{lemma}

\begin{proof}
It suffices to prove that
$\chi_{(c)}(\theta_A)=\chi_{(c)}(\theta_G)$.
 Refining our stratifications,
by Lemma \ref{refine}, we may suppose that $\mathcal{A}$ contains
at most one stratum $C/A$ with non-empty conjugation domain
$Con(A)$, and we may
 restrict $\mathcal{G}$ to a stratification of $A\times_S\A_S^m$.
Write $\mathcal{G}$ as $<\A^m_S\times_S \A^m_S,D_i/G_i,
Con(G_i)>$. Let $W$ be the support of $\mathcal{G}$, i.e. the
union of the strata $G_i$ with non-empty conjugation domain, and
  let $\pi:W\rightarrow A$ be the projection.

First, suppose that the support of $\Aa$ is empty, i.e. that
$Con(A)$ is empty. Since $\mathcal{G}$ defines the graph of a
bijection between $\Aa$ and $\mathcal{B}$, this means that
$Z(\mathcal{G},x,M)$ is empty, for any point $x$ of $S$, and any
pseudo-finite field $M$ containing $k(x)$. By Lemma \ref{decomp2},
 this implies that the support of $\mathcal{G}$ is empty. Hence,
$\chi_{(c)}(\theta_A)=\chi_{(c)}(\theta_G)=0$.

 So we may assume that
that $Con(A)$ is not empty.
  By Lemma \ref{decomp2}, the union of
the sets $Z(\mathcal{A},x,M)$, where $x$ runs over the points of
$A$, and $M$ runs over the pseudo-finite field extensions of
$k(x)$,
 is dense in $A$.
 Since $\pi$ induces, for each $x$ and each $M$, a bijection between
$Z(\mathcal{A},x,M)$ and $Z(\mathcal{G},x,M)$,
 this shows that the image of $\pi$ is dense in $A$.

Now suppose that $G$ is a stratum of $W$, and $a$ is a closed
point of $A$, such that the fiber of $G$ over $a$ has dimension
$>0$. After a refinement, we may suppose that $G$ is mapped to
$a$. Let $x$ be the image of $a$ in $S$. By Lemma \ref{decomp2},
the conjugation domain of $G$ must be empty; if not,
$Z(\mathcal{G},x,M)$ would contain infinitely many points lying
over $a$, for some pseudo-finite field $M$ containing $k(x)$.

Hence, there exists an open dense subscheme $A'$ of $A$, such that
the restriction of $\mathcal{G}$ to $\pi^{-1}(A')$ satisfies the
following property: for every stratum $G_i$ of this restriction,
$\pi:G_i\rightarrow A'$ is \'etale and finite. By Noetherian
induction, we might as well assume that $A=A'$.

We can dominate the \'etale cover $C/A$ and all the composed
covers $D_i/A$ by a common Galois cover $D/A$. By Lemma
\ref{inflate}, we might as well assume that $D=C=D_i$ for all $i$.
 We will show that
\begin{equation}\label{char}\alpha_{Con(A)}=
\sum_i Ind^{G(D/A)}_{G(D/G_i)}\alpha_{Con(G_i)}
\end{equation}
This will complete the proof, by Lemma \ref{welldef2}.

 Let $M$ be any field, and let $a$ be
an $M$-valued point on $A$. The point $a$ lifts to an $M$-valued
point $b$ on $G_i$, for some $i$, iff $C_{D/A}(a)\cap G(D/G_i)$ is
not empty, by Lemma \ref{decomp0}. In this case, there exists a
group $H$ in this intersection, such that $C_{D/G_i}(b)$ is the
conjugation class of $H$ in $G(D/G_i)$.

Conversely, if $H'$ is a subgroup of $G(D/G_i)$ for some $i$, and
if $b$ is an $M$-valued point on $G_i$ with $H'\in C_{D/G_i}(b)$,
then the image of $b$ in $A$ is an $M$-valued point $a$ whose
decomposition class $C_{D/A}(a)$ is the conjugation class of $H'$
in $G(D/A)$.

Let $g$ be any element of $G(D/A)$, denote by $H$ the subgroup of
$G(D/A)$ generated by $g$, and denote by $C_H$ its conjugation
class in $G(D/A)$.
 By definition,
$\alpha_{Con(A)}(g)=1$ if there exists a point $x$ on $S$, a
pseudo-finite field $M$ containing $k(x)$, and a point $a$ in
$Z(\Aa,x,M)$ with $H\in C_{D/A}(a)$. Else, $\alpha_{Con(A)}(g)=0$.

In the latter case, $H\notin Con(A)$, and since $\mathcal{G}$
defines a bijection between $\mathcal{A}$ and $\mathcal{B}$,
$C_H\cap Con(G_i)$ is empty, for all $i$. This means that both
members of $(\ref{char})$ vanish when evaluated in $g$.

In the first case, there exists a unique index $i$ such that $a$
lifts to an $M$-valued point $b$ on $G_i$, with
$C_{D/G_i}(b)\subset Con(G_i)$, and this point $b$ is also unique.
This means that there exists a unique index $i$ such that
$C_{D/A}(a)\cap \mathscr{P}(G(D/G_i))$ is non-empty and contained
in $Con(G_i)$, and moreover, $C_{D/A}(a)\cap Con(G_i)$ consists of
a single conjugation class in $G(D/G_i)$. Hence, both members of
$(\ref{char})$ are equal to $1$ when evaluated in $g$.
\end{proof}

By Lemma \ref{welldef} and the Elimination Theorem \ref{rel}, the
map $\chi_{(c)}$ from Definition \ref{Theta} factors through a map
of \textit{sets}
$$\chi_{(c)}:K_0(PFF_S)\rightarrow K_0^{mot}(Var_S)\otimes \Q$$
We now show that it respects the ring structures. By Lemma
\ref{welldef}, we may freely identify Galois formulas $\theta$
with equivalent Galois formulas or equivalent ring formulas. In
particular, $\chi_{(c)}$ is well-defined on $\theta_1\times
\theta_2$, $\theta_1\vee \theta_2$ and $\theta_1\wedge \theta_2$,
for any pair of Galois formulas $\theta_1$, $\theta_2$.
\begin{lemma}\label{product}
Let $\theta_1$ and $\theta_2$ be quantifier-free Galois formulas
corresponding to Galois stratifications $\mathcal{A}_1$ and
$\mathcal{A}_2$, and denote by $\theta_1\times_S \theta_2$ the
Galois formula corresponding to the product $\Aa_1\times_S\Aa_2$.
Then $$\chi_{(c)}(\theta_1\times
\theta_2)=\chi_{(c)}(\theta_1).\chi_{(c)}(\theta_2)
$$
\end{lemma}
\begin{proof}
Let $X_1,X_2$ be affine normal irreducible varieties over $S$, let
$Y_1/X_1$ and $Y_2/X_2$ be Galois covers, and let $Con_1$, resp.
$Con_2$ be conjugation domains for these covers. Put
$G_i=G(Y_i/X_i)$ for $i=1,2$.

Let $Z$ be any connected component of $Y_1\times_S Y_2$, and let
$Con$ be the conjugation domain for the cover $Z/X_1\times_S X_2$
defined in Section \ref{galoisprod}. Put $$G=G(Z/X_1\times_S
X_2)$$ and denote by $\psi:G_1\times G_2\rightarrow G$ the
projection. Observe that
$$Res_{\psi}\alpha_{Con}=\alpha_{Con_1}.\alpha_{Con_2}$$ in
$C(G_1\times G_2,\Q)$. By Lemma \ref{prod} and Lemma
\ref{welldef3}, $$M_{c,G}(Z)_{\alpha_{Con}}\cong
M_{c,G_1}(Y_1)_{\alpha_{Con_1}}\otimes
M_{c,G_2}(Y_2)_{\alpha_{Con_2}}
$$
\end{proof}

\begin{lemma}\label{sumprod}
For any pair of quantifier-free Galois formulas
$\theta_1,\,\theta_2$ with the same free variables, we have
$$\chi_{(c)}(\theta_1\vee\theta_2)+\chi_{(c)}(\theta_1\wedge
\theta_2)=\chi_{(c)}(\theta_1)+\chi_{(c)}(\theta_2)$$ For any pair
of quantifier-free Galois formulas $\theta_1,\,\theta_2$ with
disjoint sets of free variables, we have
$$\chi_{(c)}(\theta_1\wedge
\theta_2)=\chi_{(c)}(\theta_1).\chi_{(c)}(\theta_2)$$
\end{lemma}
\begin{proof}
Let $\theta_1$ and $\theta_2$ be Galois formulas corresponding to
Galois stratifications $\mathcal{A}_1$ and $\mathcal{A}_2$.

First, suppose $\theta_1$ and $\theta_2$ have the same free
variables. After a refinement, and using Lemma \ref{refine}, we
may suppose that the underlying stratifications of $\mathcal{A}_1$
and $\mathcal{A}_2$ coincide. If $Z(\theta_1\wedge \theta_2,x,M)$
is empty, for any point $x$ of $S$ and any pseudo-finite field
extension $M$ of $k(x)$, the conjugation domains of
$\mathcal{A}_1$ and $\mathcal{A}_2$ are disjoint, for each
stratum, by Corollary \ref{generic}. In this case, it is easy to
see that
$$\chi_{(c)}(\theta_1\vee\theta_2)=\chi_{(c)}(\theta_1)+\chi_{(c)}(\theta_2)$$
In general, we have
\begin{eqnarray*}
\chi_{(c)}(\theta_1\vee\theta_2)&=&\chi_{(c)}(\theta_1\wedge\theta_2)+\chi_{(c)}(\neg\theta_1\wedge\theta_2)
+\chi_{(c)}(\theta_1\wedge\neg\theta_2)\\
\chi_{(c)}(\theta_1)&=&\chi_{(c)}(\theta_1\wedge\theta_2)+\chi_{(c)}(\theta_1\wedge\neg\theta_2)
\\
\chi_{(c)}(\theta_2)&=&\chi_{(c)}(\theta_1\wedge\theta_2)+\chi_{(c)}(\neg\theta_1\wedge\theta_2)
\end{eqnarray*}
Hence,
$$\chi_{(c)}(\theta_1\vee\theta_2)+\chi_{(c)}(\theta_1\wedge
\theta_2)=\chi_{(c)}(\theta_1)+\chi_{(c)}(\theta_2)$$ if
$\theta_1$ and $\theta_2$ have the same free variables.

If $\theta_1$ and $\theta_2$ have disjoint sets of free variables,
the formula $\theta_1\wedge\theta_2$ corresponds to the product
$\mathcal{A}_1\times_S\mathcal{A}_2$. By Lemma \ref{product}, we
have
$$\chi_{(c)}(\theta_1\wedge\theta_2)=\chi_{(c)}(\theta_1).\chi_{(c)}(\theta_2)$$
\end{proof}

In analogy with \cite{DL1}, Theorem 2.1, we state
\begin{theorem}\label{morph}
There exists a unique ring morphism
$$\chi_{(c)}:K_0(PFF_S)\rightarrow K_0^{mot}(Var_S)\otimes \Q$$
satisfying the following two properties:
\\(i) for every quantifier-free ring formula $\varphi$, the image
of $[\varphi]$ under $\chi_{(c)}$ equals $\chi_{c}([X])$, where
$X$ is the $S$-constructible set defined by $\varphi$.
\\(ii) let $X$ be a normal affine irreducible variety over $S$,
let $Y$ be a Galois cover of $X$, and let $C$ be a cyclic subgroup
of the Galois group $G$ of the cover $Y/X$. To these data, we can
associate a ring formula $\varphi_{Y/X,C}$ over $S$ (see Section
\ref{galoisring}), whose interpretation, in any point $x$ of $S$,
and for any field $K$ containing $k(x)$, is the set of $K$-valued
points $a\in X_x(K)$ with $C\in C_{Y/X}(a)$. Then
$$\chi_{c}([\varphi_{Y/X,C}])=\frac{|C|}{|N_G(C)|}\chi_{(c)}([\varphi_{Y/(Y/C),C}])\,, $$
where $N_G(C)$ denotes the normalizer of $C$ in $G$.
\end{theorem}

\begin{proof}

 \textit{1.
Uniqueness:} By Corollary \ref{generators}, the morphism
$\chi_{(c)}$ is determined by the images of classes of the form
$[\varphi_{Y/X,C}]$, with $C$ cyclic. Hence, by $(ii)$, $\chi_{c}$
is determined by the images of classes of the form
$[\varphi_{Y/(Y/C),C}]$, where $C$ is cyclic. However, by (i) and
(ii),
$$|C|\,\chi_{c}([Y/C])=\sum_{A\ subgroup\ of\ C}|A|\,\chi_{(c)}([\varphi_{Y/(Y/A),A}])$$
since the formulas $\varphi_{Y/(Y/C),A}$ yield a partition of
$Y/C$, and since
$$|C|\,\chi_{(c)}([\varphi_{Y/(Y/C),A}])= |A|\,\chi_{(c)}([\varphi_{Y/(Y/A),A}])$$
by property (ii). This recursion formula determines
$\chi_{(c)}([\varphi_{Y/(Y/C),C}])$.

\textit{2. Existence:} By Lemma \ref{sumprod} and Lemma
\ref{welldef}, the map $\chi_{(c)}$ from definition \ref{Theta}
factors to a ring morphism
$$\chi_{(c)}:K_0(PFF_S)\rightarrow K_0^{mot}(Var_S)\otimes
\Q$$
 The fact that property
(i) is satisfied, follows immediately from the definition. Let us
prove property (ii). Put $G=G(Y/X)$, and denote by
$\psi:C\rightarrow G$ the inclusion. We have $Y_C=Res_{\psi}Y_G$,
where $Y_C$ is the variety $Y$ with good $C$-action as a Galois
cover of $Y/C$, and $Y_G$ is the variety $Y$ with good $G$-action
as a Galois cover of $X$. If $Con(Y/X)$ is the conjugation class
of $C$ in $G$, and if we denote by $Con(Y/(Y/C))$ the conjugation
domain $\{C\}$, then
$$\alpha_{Con(Y/X)}=\frac{|C|}{|N_G(C)|}Ind_{\psi}\alpha_{Con(Y/(Y/C))}$$
Hence, we may conclude by Lemma \ref{welldef2}.
\end{proof}
\begin{cor}
There exists a canonical ring morphism
$$K_0(PFF_S)\rightarrow K_0(D^b(S,\Q_\ell))\otimes \Q$$
mapping a ring formula $\varphi$ defining a locally closed subset
$X$ of $\A^m_S$, to the class of $R\pi_!(\Q_\ell)$, where
$\pi:\A^m_S\rightarrow S$ is the structural morphism.
\end{cor}
\begin{proof}
The morphism is obtained by composing $\chi_{(c)}$ with the
realization morphism from Corollary \ref{realiz}. By point (i) of
Theorem \ref{morph}, the image of $[\varphi]$ is nothing but the
image of $\chi_{(c)}(X)$ in $K_0(D^b(S,\Q_\ell))$. Let us show
that this image coincides with the class of $R\pi_!(\Q_\ell)$. By
additivity of $R\pi_!$, and constructible resolution of
singularities, we may suppose that $X$ is smooth and proper over
$S$. In this case, the result follows simply by definition of the
realization functor $\mathcal{M}^o_+(S)\rightarrow
D^b(S,\Q_\ell)$.
\end{proof}

\begin{prop}[Base Change]\label{basechange}
Let $T$ be another $k$-variety, with a morphism $f:T\rightarrow
S$. This morphism induces base change morphisms
$K_0(PFF_S)\rightarrow K_0(PFF_T)$ and $K^{mot}_0(Var_S)\otimes
\Q\rightarrow K_0^{mot}(Var_T)\otimes \Q$, which we both denote by
$f^{*}$. We get a commutative diagram
$$\begin{CD}
K_0(PFF_S)@>\chi_{(c)}>>K_0^{mot}(Var_S)\otimes \Q
\\@Vf^{*}VV @Vf^{*}VV
\\K_0(PFF_T)@>\chi_{(c)}>>K_0^{mot}(Var_T)\otimes \Q\,.
\end{CD}$$
\end{prop}
\begin{proof}
Let $\varphi(x_1,\ldots,x_m)$ be a ring formula over $S$, and let
$\varphi'$ be its pullback to $T$. Suppose that $\varphi$ is
equivalent to a quantifier-free Galois formula $\theta$, with
corresponding Galois stratification
$\mathcal{A}=<A_i,C_i/A_i,Con(A_i)>$, i.e.
$Z(\varphi,x,M)=Z(\mathcal{A},x,M)$ for each point $x$ of $S$, and
each pseudo-finite field extension $M$ of the residue field
$k(x)$. Let $\theta'$ be the Galois formula, associated to a
Galois stratification of $\A^{m}_T$,
 obtained by pulling back $\mathcal{A}$ via the morphism
$\A^{m}_T\rightarrow \A^{m}_S$ induced by $f$, as explained in
Section \ref{galoispback}. We denote this Galois stratification by
$\mathcal{A}'=<A'_j,C'_j/A'_j, Con(A'_j)>$. It is clear that
$\varphi'$ is equivalent to $\theta'$.

 Now we compare $\chi_{(c)}(\theta)$ and $\chi_{(c)}(\theta')$. We may
suppose that the support of $\mathcal{A}$ consists of a single
stratum $A$, with Galois cover $C/A$ and conjugation domain
$Con(A)$. By a stratification argument, we may assume that $T$ is
irreducible and smooth, and that $A':=A\times_k T$ is a stratum of
$\mathcal{A}'$. It suffices to prove that
$$f^{*}\chi_{c,\alpha_{Con(A)}}([C])
=\chi_{c,\alpha_{Con(A')}}([C'])$$ where $C'/A'$ and $Con(A')$ are
the induced Galois cover and conjugation domain; i.e. $C'$ is a
connected component of $C\times_A A'$, and $Con(A')$ is the set of
elements of $Con(A)$ which are contained in $G(C'/A')$. We denote
by $\psi:G(C'/A')\rightarrow G(C/A)$ the inclusion. By Lemma
\ref{basechange0},
$$f^{*}\chi_{c,\alpha_{Con(A)}}([C])=
\chi_{c,\alpha_{Con(A)}}([C\times_{S}T])$$ while $C\times_S
T=Ind_{\psi}C'$ and $\alpha_{Con(A')}=Res_{\psi}\alpha_{Con(A)}$.
So we can conclude by applying Lemma \ref{welldef3}.
\end{proof}

\section{The relative motivic Poincar\'e series}\label{series}
Let $S$ be a variety over a field $k$ of characteristic zero. Let
$X$ be a separated scheme, of finite type over $S$. For any
integer $n> 0$, we denote by $S[t]/(t^{n})$ the $S$-scheme
$S\times_k \mathrm{Spec}\,k[t]/(t^n)$. We define, for each integer
$n\geq 0$, a contravariant functor
$$F_n:(Sch_S)^{op}\rightarrow (Sets)$$ from the category $(Sch_S)$ of $S$-schemes, to the category
of sets, by
$$F_n(Y)=Hom_{Sch_S}(Y\times_{S}S[t]/(t^{n+1}),X)$$
For any field $K$ containing $k$, we call the points of
$F_n(\mathrm{Spec}\,K)$ the $K$-valued relative $n$-jets on $X/S$.
\begin{prop}
The functor $F_n$ is representable by a separated $S$-scheme of
finite type $\mathcal{L}_n(X/S)$, for each $n$.
\end{prop}
\begin{proof}
The proof is analogous to the proof in the absolute case
$S=\mathrm{Spec}\,k$. Suppose that $X$ is affine over $S$, with
 $X=\mathrm{Spec}\,\mathcal{O}_S[x_1,\ldots,x_m]/(f_1,\ldots,f_r)$.
Let $$a=(a_{1,0}+a_{1,1}t+\ldots+a_{1,n}t^{n},\ldots,a_{m,0}+
\ldots,a_{m,n}t^{n})$$ be an $m$-tuple of elements of
$\mathcal{O}_S[t]/(t^{n+1})$. The system of equations
$f_i(a)\equiv 0\,\mod\,t^{n+1}$, for $i=1,\ldots,r$, puts
algebraic conditions on the coefficients $a_{i,j}$, and if we
consider these coefficients as affine coordinates on
$\A^{m(n+1)}_S$, the $f_i$ define a closed subscheme
$\mathcal{L}_n(X/S)$ of $\A^{m(n+1)}_S$. A gluing procedure yields
the general case.
\end{proof}

Note that $\mathcal{L}_0(X/S)\cong X$. By Yoneda's Lemma, the
truncation map $S[t]/(t^{n+1})\rightarrow S[t]/(t^n)$ induces a
truncation morphism of $S$-schemes
$$\pi^{n+1}_n:\mathcal{L}_{n+1}(X/S)\rightarrow
\mathcal{L}_n(X/S)$$ for each $n\geq 0$. These morphisms are
affine, and hence, we can take the projective limit in the
category of $S$-schemes to obtain an $S$-scheme
$\mathcal{L}(X/S)$.
 It comes with natural
projections $\pi_{n}:\mathcal{L}(X/S)\rightarrow
\mathcal{L}_{n}(X/S)$. For any field $K$ containing $k$, we call
the points of $\mathcal{L}(X/S)(K)$ the $K$-valued relative arcs
on $X/S$.
%

 When $X$ is smooth over
$S$, the morphisms $\pi_{n}^{n+1}$ are piecewisely trivial
fibrations with fiber $\A^{d}_{S}$, where $d$ is the relative
dimension of $X$. An $S$-morphism $h$ from $X$ to $X'$ induces a
morphism $h$ from $\mathcal{L}(X/S)$ to $\mathcal{L}(X'/S)$ by
composition, and a morphism of $k$-varieties $f:W\rightarrow S$
induces a pull-back morphism $f^{*}:\mathcal{L}(X/S)\rightarrow
\mathcal{L}(X\times_SW/W)$. If $Y$ is a separated $W$-scheme of
finite type, this morphism $f$ induces also a forgetful morphism
$f_{*}:\mathcal{L}(Y/W)\rightarrow \mathcal{L}(Y/S)$.

\begin{lemma}\label{bchange}
For any $k$-variety $T$, endowed with a morphism $f:T\rightarrow
S$, and for any separated $S$-scheme $X$ of finite type, there are
canonical isomorphisms of $T$-schemes
\begin{eqnarray*}
\mathcal{L}_n(X\times_S T/T)&\cong&\mathcal{L}_n(X/S)\times_S T,
\\\mathcal{L}(X\times_S T/T)&\cong&\mathcal{L}(X/S)\times_S T,
\end{eqnarray*}
compatible with the truncation morphisms $\pi^{m}_{n}$ and
$\pi_n$.
\end{lemma}
\begin{proof}
Denote, for any variety $U$, the category of $U$-schemes by
$Sch_U$.
 The scheme $\mathcal{L}_n(X/S)\times_S T$ represent the
 functor $$G_n: (Sch_T)^{op}\rightarrow (Sets)$$
 mapping a $T$-scheme
$Y$ to \begin{eqnarray*} Hom_{Sch_T}(Y,\mathcal{L}_n(X/S)\times_S
T)&=&Hom_{Sch_S}(Y,\mathcal{L}_n(X/S)) \\ &=&
Hom_{Sch_S}(Y\times_S S[t]/(t^{n+1}),X).\end{eqnarray*} There is a
canonical isomorphism of $T$-schemes $$Y\times_S
S[t]/(t^{n+1})\cong Y\times_T T[t]/(t^{n+1})$$ and we have
$$Hom_{Sch_S}(Y\times_S S[t]/(t^{n+1}),X)=Hom_{Sch_T}(Y\times_T T[t]/(t^{n+1}),X\times_S
T)$$ Hence, there is a natural equivalence of functors between
$G_n$ and the functor represented by $\mathcal{L}_n(X\times_S
T/T)$, and, by Yoneda's lemma, it corresponds to an isomorphism
$$\mathcal{L}_n(X\times_S T/T)\cong\mathcal{L}_n(X/S)\times_S T$$
Taking projective limits over $n$ concludes the proof.
\end{proof}

For any scheme $Y$, we denote by $Y_{red}$ the underlying reduced
scheme.
\begin{lemma}\label{reduced}
For any separated $S$-scheme $X$ of finite type, there exists an
integer $a>0$, such that, for each $n\geq 0$, the truncation
morphism
$$\mathcal{L}_{na}(X/S)_{red}\rightarrow \mathcal{L}_{n}(X/S)$$ factors
through a morphism
$$\mathcal{L}_{na}(X/S)_{red}\rightarrow
\mathcal{L}_{n}(X_{red}/S)$$
\end{lemma}
\begin{proof}
We may assume that $X$ and $Y:=\mathcal{L}_{na}(X/S)_{red}$ are
affine, say $X=\mathrm{Spec}\,B$, and $Y=\mathrm{Spec}\,A$. Denote
by $\mathcal{N}_B$ the nilradical of $B$. We have to prove that
there exists an integer $a>0$, such that, for each $n\geq 0$, for
each reduced $k$-algebra $A$, and for each morphism of
$k$-algebras
$$g:B\rightarrow  A[t]/(t^{na+1})$$ the composition
$h:B\rightarrow  A[t]/(t^{n+1})$ factors through
$B/\mathcal{N}_B$.

In other words, we have to show that the ideal $g^{-1}(t^{n+1})$
contains $\mathcal{N}_B$. Since $B$ is Noetherian, we can choose a
positive integer $a$, such that $b^{a}=0$ for any element $b$ of
$\mathcal{N}_B$. Since $A$ is reduced, $g(b)$ has to belong to the
ideal $(t^{n+1})$, as soon as $b\in\mathcal{N}_B$.
\end{proof}

For any point $x$ on $S$ and any separated $S$-scheme of finite
type $X$, the fiber of $\mathcal{L}(X/S)$ over $x$ is canonically
isomorphic to $\mathcal{L}(X\times_S x)$, where $\mathcal{L}(.)$
denotes the arc scheme as defined in \cite{DLinvent},\,p.1. Hence,
for any field $K$, giving a $K$-valued relative arc $\psi$ on
$X/S$, amounts to giving an arc in $(X\times_x S)(K[[t]])$, where
$x$ denotes the image of $\pi_0(\psi)$ in $S$.

For any separated $S$-scheme $Y$ of finite type, we put
$[Y]:=[Y_{red}]$ in $K_0(Var_S)$.

\begin{definition}
If $X$ is a separated $S$-scheme of finite type, we define its
relative Igusa Poincar\'e series as
$$Q(X/S;T)=\sum_{n=0}^{\infty}[\mathcal{L}_n(X/S)]T^n\qquad\in K_0(Var_S)[[T]]$$
\end{definition}

 We can also generalize Denef and Loeser's definition
of the geometric and arithmetic Poincar\'e series
\cite{DL5,DL,DL3,DL1,DL2} to our relative setting. The aim is to
define $P_{geom}(X/S;T)$ and $P_{arith}(X/S;T)$ as objects in
$K_0(Var_S)[[T]]$, resp.
\\ $(K_0^{mot}(Var_S)\otimes \Q)[[T]]$, such that base-change to any
point $x$ of $S$, yields the ``classical'' motivic Poincar\'e
series of the fiber $X\times_S x$.

First, we need to prove a uniform version of Greenberg's Theorem
\cite{Gr}. We will show that there exist positive integers $c,e$,
such that, for any point $x$ on $S$, for any positive integer $n$,
and for any field extension $K$ of $k(x)$, the projections $\pi_n
\mathcal{L}(X_x)(K)$ and $\pi_n^{cn+e}\mathcal{L}_{cn+e}(X_x)(K)$
coincide, where $X_x$ denotes the fiber $X\times_S x$ of $X$ over
$x$. In other words, a $K$-valued $n$-jet on a fiber $X_x$ lifts
to a $K$-valued arc on $X_x$, as soon as it lifts to a $K$-valued
$(cn+e)$-jet on $X_x$. Using our relative arc and jet spaces, we
can write this property as
$$\pi_n(\mathcal{L}(X/S)(K))=\pi^{cn+e}_n(\mathcal{L}_{cn+e}(X/S)(K))$$

We will proceed in two steps. We start by giving a short
alternative proof of the absolute case, using resolution of
singularities. Then, we show how to extend this proof to obtain
the uniform version stated above. First of all, we recall the
following definition.

\begin{definition}
Let $X$ be a separated scheme of finite type over $k$, and let
$\mathcal{I}$ be an ideal sheaf on $X$. For any field $K$
containing $x$, and any arc $\psi$ in $\mathcal{L}(X)(K)$, we put
$x:=\pi_0(\psi)\in X$, and we define the order of $\mathcal{I}$ at
$\psi$ by
$$ord_{\mathcal{I}}(\psi):=\min \{ord_t f(\psi)\,|\,f\in \mathcal{I}_x \}\in \N\cup\{\infty\}$$
If $Z$ is a closed subscheme of $X$, we denote the defining ideal
sheaf of $Z$ in $X$ by $\mathcal{I}_Z$, and we put
$$ord_Z(\psi):=ord_{\mathcal{I}_Z}(\psi)$$
We call this value the contact order of $\psi$ with $Z$. Note that
$ord_Z(\psi)=\infty$ iff $\psi$ is contained in $Z$.
\end{definition}

\begin{theorem}[Greenberg]
Let $Z$ be a separated scheme of finite type over $k$. There exist
positive integers $c,e$, such that,
 for any positive integer $n$, and for any field extension $K$ of $k$,
the projections $\pi_n \mathcal{L}(Z)(K)$ and
$\pi_n^{cn+e}\mathcal{L}_{cn+e}(Z)(K)$ coincide.
\end{theorem}
\begin{proof}
%
We embed $Z$ as a subscheme in some smooth ambient $k$-variety
$Y$, and we take an embedded principalization $h:Y'\rightarrow Y$
for the defining ideal sheaf $\mathcal{I}_Z$ of $Z$ in the
structure sheaf $\mathcal{O}_{Y}$, defined over $k$, such that $h$
is an isomorphism over the complement $Y\setminus Z$, see
\cite[2.5]{VMayor}. This means that $h$ is a proper birational
morphism from a smooth variety $Y'$ to $Y$, such that
$$\mathcal{I}_Z.\mathcal{O}_{Y'}=\mathcal{O}_{Y'}(-\sum_{i=1}^{N}r_i E_i)$$
where $E=\sum_{i=1}^{N}E_i$ is a normal crossing divisor on $Y'$,
and where the $r_i$ are positive integers. We denote the maximum
of the multiplicities $r_i$ by $R$.

Choose an integer $n>0$. Let $K$ be a field extension of $k$, and
let $\psi$ be a $K$-valued arc on $Y$, such that $$nNR\leq
ord_{Z}(\psi)<\infty$$  This means that
$$\pi_{nNR-1}(\psi):\mathrm{Spec}\,K[t]/(t^{nNR})\rightarrow Y$$ factors through a
$(nNR-1)$-jet on $Z$. Since $h$ is an isomorphism over $Y\setminus
Z$, it follows from the valuative criteria for properness and
separateness, that there exists a unique $K$-valued arc $\psi'$ on
$Y'$, such that $h\circ \psi'=\psi$. We put
$\nu_i=ord_{E_i}(\psi')$, and we denote by $\nu$ the maximum of
the $\nu_i$. We may assume that $\nu=\nu_1$.

Now, observe that $nNR\leq ord_Z(\psi)=\sum_{i=1}^{N}r_i\nu_i\leq
\nu NR$, hence $\nu\geq n$. Since the divisor $E_1$ is smooth,
there exists an arc $\zeta'$ in $\mathcal{L}(E_1)(K)$ such that
$\pi_{\nu-1}(\zeta')=\pi_{\nu-1}(\psi')$ in
$\mathcal{L}_{\nu_1}(E_1)(K)$. If we denote by $\zeta$ the image
of $\zeta'$ under $h$, then
$\pi_{\nu-1}(\zeta)=\pi_{\nu-1}(\psi)=h(\pi_{\nu-1}(\psi'))$ in
$\mathcal{L}_{\nu-1}(Z)$, and, a fortiori,
$\pi_{n-1}(\zeta)=\pi_{n-1}(\psi)$.

Hence, we see that, for any $K$-valued $n$-jet $\psi_n$ on $Z$,
this jet $\psi_n$ lifts to a $K$-valued arc on $Z$, iff it lifts
to a $K$-valued $(nNR+NR-1)$-jet on $Z$.
\end{proof}

\begin{theorem}[Uniform version]\label{uniform}
Let $S$ be a variety over $k$, and let $X$ be a separated scheme
of finite type over $S$. There exist positive integers $c,e$, such
that, for any point $x$ on $S$, for any positive integer $n$, and
for any field extension $K$ of $k(x)$, the projections $\pi_n
\mathcal{L}(X_x)(K)$ and $\pi_n^{cn+e}\mathcal{L}_{cn+e}(X_x)(K)$
coincide, where $X_x$ denotes the fiber of $X$ over $x$.
\end{theorem}
\begin{proof}
By the proof of the previous theorem, it suffices to find, for
each fiber $X_x$, an embedding in some smooth ambient space, and
an embedded principalization for its defining ideal sheaf, such
that the numbers $N,R$ are uniformly bounded. This follows from
the existence of a constructible embedded resolution, by
Proposition \ref{resolution}.
\end{proof}

\begin{definition}
If $X$ is a separated $S$-scheme of finite type, we define its
relative geometric Poincar\'e series as
$$P_{geom}(X/S;T)=\sum_{n=0}^{\infty}[\pi_n\mathcal{L}(X/S)]T^n\qquad\in K_0(Var_S)[[T]]$$
\end{definition}

By Theorem \ref{uniform},
$(\pi_n\mathcal{L}(X/S))_{red}=(\pi^{n'}_n\mathcal{L}_{n'}(X/S))_{red}$
for some integer $n'>n$, hence the truncation
$(\pi_n\mathcal{L}(X/S))_{red}$ is, by Chevalley's Theorem, a
constructible subset of $\mathcal{L}_n(X/S)_{red}$, and its
Grothendieck bracket is well-defined.

\begin{prop}[Base Change]
Let $X$ be an $S$-variety, and let $R(X/S;T)$ be either its
relative Igusa Poincar\'e series, or its relative geometric
Poincar\'e series. Let $W\rightarrow S$ be a morphism of
$k$-varieties.
 If we take the image of the
coefficients of $R$ under the base change morphism
$K_0(Var_S)\rightarrow K_0(Var_{W})$, then we obtain the
corresponding Poincar\'e series $R(X\times_S W;T)$.
\end{prop}
\begin{proof}
This follows immediately from Lemma \ref{bchange}.
\end{proof}
In particular, by base change to any point $x$ of $S$, we recover
the absolute Poincar\'e series of the (not necessarily reduced)
fiber $X_x$.

Now, we define the relative arithmetic Poincar\'e series. If $X$
is affine over $S$, say $X$ is a closed subscheme of $\A_S^m$,
Theorem \ref{uniform} implies that we can find, for each positive
integer $n$,
 a ring formula $\psi_n(x_1,\ldots,x_{(n+1)m})$ over $S$,
such that, for any point $x$ on $S$, and any field extension $K$
of $k(x)$, the $(n+1)m$-tuples over $K$ satisfying $\psi_n$,
correspond to the points of $\pi_n (\mathcal{L}(X_x)(K))$.
 Here we identified $\pi_n (\mathcal{L}(X_x)(K))$ with a subset of
$\mathcal{L}_{n}(\A_{k(x)}^m)(K)=\A^{(n+1)m}_{k(x)}(K)$.

\begin{definition}
 We define the relative arithmetic Poincar\'e series as
$$P_{arith}(X/S;T)=\sum_{n=0}^{\infty}
\chi_{(c)}([\psi_n])T^n\qquad \in (K_0^{mot}(Var_S)\otimes
\Q)[[T]]$$
\end{definition}
If $X$ is not affine over $S$, we can still define its arithmetic
Poincar\'e series, using definable subassignements as in
\cite{DL1}.

\begin{prop}[Base Change]
Let $X$ be a separated $S$-scheme of finite type, and let
$W\rightarrow S$ be a morphism of $k$-varieties.
 If we take the image of the
coefficients of $P_{arith}(X/S;T)$ under the base change morphism
$K_0^{mot}(Var_S)\otimes \Q \rightarrow K_0^{mot}(Var_W)\otimes
\Q$, then we obtain the arithmetic Poincar\'e series
$P_{arith}(X\times_S W;T)$.
\end{prop}
\begin{proof}
This follows from Lemma \ref{basechange}.
\end{proof}

\section*{Acknowledgements}
The author would like to thank Fran\c{c}ois Loeser for proposing
the subject of this article.
\bibliographystyle{hplain}
\bibliography{wanbib,wanbib2}
\end{document}